\documentclass[twoside,12pt]{article}
\usepackage[]{amsmath,amsthm,amssymb}
\usepackage[latin1]{inputenc}

\begin{document}

\title{{\Large\bf  Generalized product formulas for Whittaker's functions and  a novel  class of index transforms}}

\author{Semyon  YAKUBOVICH}
\maketitle

\markboth{\rm \centerline{ Semyon   YAKUBOVICH}}{}
\markright{\rm \centerline{PRODUCT FORMULAS  AND INDEX   TRANSFORMS FOR WHITTAKER'S FUNCTIONS}}

\begin{abstract} {\noindent Generalized product formulas and index transforms, involving  products  of Whittaker's functions of different   indices  are established  and investigated. The corresponding inversion formulas are found. Particular cases cover index transforms  with products of the modified Bessel and Whittaker's functions. For our goals the Kontorovich-Lebedev and Olevskii transforms of a  complex index with nonzero real part are involved. }

\end{abstract}
\vspace{4mm}

{\bf Keywords}: {\it Whittaker function, modified Bessel function,  Kontorovich-Lebedev transform, Olevskii transform, Laplace transform, Riemann-Liouville fractional integrals, Mellin transform}

{\bf AMS subject classification}:  45A05,  44A15,   33C15, 33C10

\vspace{4mm}

\section {Product formulas for Whittaker's functions}

The main goal of this paper is to investigate a class of index transformations with respect to the second pure imaginary index $i\tau,\ \tau \in \mathbb{R}$ of the product of Whittaker's functions $W_{\alpha,i\tau}(x) W_{\beta,i\tau}(x)$ [2, Vol. III] of the positive argument $ x >0$ and   real indices $\alpha, \beta \in \mathbb{R}$.  Historically, index transforms, involving Whittaker's functions (index Whittaker transforms) appear for the first time in 1964 when Wimp [5] discovered the following  reciprocal pair of integral transformations

$$F(\tau)=  \int_0^\infty   W_{\mu,i\tau}(x)  f(x) {dx\over x^2},\  \tau >0,\eqno(1.1)$$

$$f(x)= { 1 \over \pi^2  } \int_0^\infty  \tau \sinh(2\pi\tau) \left|\Gamma\left({1\over 2}- \mu + i\tau\right)\right|^2  W_{\mu,i\tau}(x)  F(\tau) d\tau,\ x >0,\eqno(1.2)$$
 where $\mu \in \mathbb{R}$ and $\Gamma(z)$ is the Euler gamma function [2, Vol. III].    As we see in (1.1), (1.2), the corresponding integrals depend upon a second parameter (index) of the Whittaker function and involve it as the variable of integration.  This function generalizes the modified Bessel or Macdonald function $K_\nu(x)$ [2, Vol. II] via the equality 
 
 $$W_{0,i\tau}(2x) = \sqrt {{2x\over \pi}}\ K_{i\tau}(x),\eqno(1.3)$$ 
and  the pair (1.1), (1.2) reduces to Kontorovich-Lebedev transforms [6].  Generally,  these functions are related through the following  integral (cf. Entry 2.16.8.4 in [2,  Vol. II])

 $$ \Gamma\left({1\over 2}- \mu + i\tau\right)\Gamma\left({1\over 2}- \mu - i\tau\right) W_{\mu,i\tau}(x) $$
 
 $$= 2(4x)^\mu e^{-x/2} \int_0^\infty  t^{-2\mu} e^{-t^2/(4x)} K_{2i\tau}(t) dt,\ x >0, \  {\rm Re} (\mu) < {1\over 2},\ \tau \in \mathbb{R}.\eqno(1.4)$$
The Whittaker function $W_{\alpha,\nu} (x )$  is, by definition, the solution of Whittaker's  differential equation

$${d^2 u\over dx^2} + \left({\alpha\over x} - {1\over 4} + {1/4-\nu^2\over x^2} \right) u = 0.\eqno(1.5)$$
Concerning the product of the Whittaker functions, recently we established the product formula for the case $\alpha=\beta,\ \mu=\nu$ and different arguments (see [3, formula (3.19)])

$$W_{\alpha, \nu}(x) W_{\alpha, \nu}(y) = \int_0^\infty   W_{\alpha, \nu}(\xi) k_\alpha (x,y,\xi) {d\xi\over \xi^2},\quad x,y >0,\ \alpha,\nu \in \mathbb{C},\eqno(1.6)$$
where

$$k_\alpha (x,y,\xi) = 2^{-1-\alpha} \pi^{-1/2} (xy\xi)^{1/2} \exp\left( {x\over 2}+ {y\over 2}+ {\xi\over 2} - { (xy+x\xi+y\xi)^2\over 8xy\xi} \right)$$

$$\times  D_{2\alpha} \left(  { xy+x\xi+y\xi\over (2xy\xi)^{1/2}}\right)$$
and $D_\mu(z)$ is the parabolic cylinder function  [1, Vol. II].      It gives a departure point for the corresponding convolution theory for the index Whittaker transforms.  However, for our further purposes in the sequel we will generalize the product formula (3.26) in [3] in terms of the Gauss hypergeometric function  [2, Vol. III]

$$W_{\alpha, i\tau}(x) W_{\alpha, i\tau}(y) = {(xy)^\alpha\ e^{-(x+y)/2}\over \Gamma(1-2\alpha)} \int_0^\infty e^{-t} t^{-2\alpha} $$

$$\times  {}_2F_1 \left( {1\over 2} -\alpha-i\tau,\  {1\over 2} -\alpha+ i\tau ; 1-2\alpha ; - \left({1\over x}+ {1\over y} \right) t - {t^2\over xy} \right) dt,\eqno(1.7)$$
where $x,y >0,\ \tau \in \mathbb{R},\ \alpha < 1/2$.
To do this, we will employ the expression of the Whittaker function in terms of the Tricomi function or the Kummer confluent hypergeometric function of the second kind  [2, Vol. III]

$$\Psi(a,b; x) = e^{x/2} x^{-b/2}\  W_{b/2-a, (b-1)/2} (x).$$
Indeed, combining with Entry 2.21.2.17 in [2, Vol. III], we derive immediately the following product formula

$$ W_{\alpha,i\tau}(x) W_{\beta,i\tau}(y) = (xy)^{1/2+i\tau}  \ e^{-(x+y)/2}\ \Psi\left({1\over 2} -\alpha +i\tau, 1+2i\tau; x \right)$$

$$\times  \ \Psi\left({1\over 2} -\beta +i\tau, 1+2i\tau; y \right) = { (xy)^{1/2+i\tau}  \ e^{-(x+y)/2} \over \Gamma(1-\alpha-\beta)}  \int_0^\infty e^{-t} t^{-\alpha-\beta} $$

$$\times (t+x)^{\alpha-1/2-i\tau} (t+y)^{\beta-1/2-i\tau}$$

$$\times   {}_2F_1 \left( {1\over 2} -\alpha+i\tau,\  {1\over 2} -\beta+ i\tau ; 1-\alpha-\beta ; 1- {xy\over (t+x)(t+y)}  \right) dt,\quad \alpha+\beta < 1.\eqno(1.8)$$
Meanwhile, Entry 7.3.1.3 in [2, Vol. III] for the Gauss hypergeometric function suggests the following equalities from (1.8)

$$  W_{\alpha,i\tau}(x) W_{\beta,i\tau}(y) =  { (xy)^{\alpha}  \ e^{-(x+y)/2} \over \Gamma(1-\alpha-\beta)}  \int_0^\infty e^{-t} t^{-\alpha-\beta}  (t+y)^{\beta-\alpha} $$

$$\times  {}_2F_1 \left( {1\over 2} -\alpha+i\tau,\  {1\over 2} -\alpha- i\tau ; 1-\alpha-\beta ; - t \left({1\over x}+ {1\over y}\right)- {t^2\over xy} \right)  dt$$

$$=  { (xy)^{\beta}  \ e^{-(x+y)/2} \over \Gamma(1-\alpha-\beta)}  \int_0^\infty e^{-t} t^{-\alpha-\beta}  (t+x)^{\alpha-\beta}$$

$$\times    {}_2F_1 \left( {1\over 2} -\beta+i\tau,\  {1\over 2} -\beta- i\tau ; 1-\alpha-\beta ; - t \left({1\over x}+ {1\over y}\right)- {t^2\over xy} \right)  dt.\eqno(1.9)$$
In particular, for $x=y$ after a simple substitution it gives the following Laplace integrals for the Gauss hypergeometric functions

$$  W_{\alpha,i\tau}(x) W_{\beta,i\tau}(x) =  { x  \ e^{- x} \over \Gamma(1-\alpha-\beta)}  \int_0^\infty e^{- xt} t^{-\alpha-\beta}  (t+1)^{\beta-\alpha} $$

$$\times  {}_2F_1 \left( {1\over 2} -\alpha+i\tau,\  {1\over 2} -\alpha- i\tau ; 1-\alpha-\beta ; -  2t -  t^2 \right)  dt$$

$$=  { x \ e^{- x} \over \Gamma(1-\alpha-\beta)}  \int_0^\infty e^{- xt} t^{-\alpha-\beta}  (t+1)^{\alpha-\beta}$$

$$\times    {}_2F_1 \left( {1\over 2} -\beta+i\tau,\  {1\over 2} -\beta- i\tau ; 1-\alpha-\beta ;  - 2t - t^2 \right)  dt.\eqno(1.10)$$

In order to extend the product formula for Whittaker functions with real  parameters $\alpha, \beta$ and two complex indices $\mu-i\tau,\ \nu+i\tau,\ \mu,\nu, \tau \in \mathbb{R}$ we employ Entry 2.16.8.4 in [2, Vol. II] to write

$$ \Gamma\left({1\over 2}-\alpha+\mu-i\tau\right)  \Gamma\left({1\over 2}-\alpha-\mu+i\tau\right)   \Gamma\left({1\over 2}-\beta+\nu+i\tau\right)  \Gamma\left({1\over 2}-\beta-\nu-i\tau\right) $$

$$\times x^{\alpha+\beta}  e^{1/x} W_{\alpha, \mu-i\tau}\left({1\over x}\right) W_{\beta, \nu+i\tau}\left({1\over x}\right) $$

$$= 4 \int_0^\infty  \int_0^\infty e^{-x(s^2+t^2)} s^{-2\alpha} t^{- 2\beta}  K_{2(\mu-i\tau)} \left(  s \right) K_{2(\nu+i\tau)} \left(  t \right)dsdt\eqno(1.11)$$
under the convergence conditions $x >0, \ \tau  \in \mathbb{R}, ,\ \alpha, \beta < 1/2, \  1-2\alpha >  2| \mu|,\  1-2\beta >  2|\nu|$. Meanwhile,  the product of the Macdonald functions can be expressed via relation 2.16.13.4 in [2, Vol. II], and we have 

$$2 K_{2(\mu-i\tau)} \left(  s \right) K_{2(\nu+i\tau)} \left(  t \right) = \int_{-\infty}^\infty  e^{-2i\tau y+ (\mu-\nu)y}  \left({s+te^y\over s e^y+t}\right)^{\mu+\nu}$$

$$\times K_{2(\mu+\nu)} \left( \sqrt{ s^2+t^2 + 2st \cosh y}\right) dy,\quad s,t > 0.\eqno(1.12)$$
Plugging the right-hand side of the equality (1.12) into (1.11),  we interchange the order of integration in the obtained iterated integral via the dominated convergence theorem because 

$$\int_0^\infty  \int_0^\infty e^{-x(s^2+t^2)} s^{-2\alpha} t^{- 2\beta}   \int_{-\infty}^\infty  e^{(\mu-\nu)y}  \left({s+te^y\over s e^y+t}\right)^{\mu+\nu}$$

$$\times K_{2(\mu+\nu)} \left( \sqrt{ s^2+t^2 + 2st \cosh y}\right) dy ds dt$$

$$= 2 \int_0^\infty  \int_0^\infty e^{-x(s^2+t^2)} s^{-2\alpha} t^{- 2\beta} K_{2\mu} \left(  s \right) K_{2\nu} \left(  t \right) ds dt < \infty,\ x > 0.$$
Therefore, using polar coordinates, equality (1.11) becomes

$$ \Gamma\left({1\over 2}-\alpha+\mu-i\tau\right)  \Gamma\left({1\over 2}-\alpha-\mu+i\tau\right)   \Gamma\left({1\over 2}-\beta+\nu+i\tau\right)  \Gamma\left({1\over 2}-\beta-\nu-i\tau\right) $$

$$\times x^{\alpha+\beta}  e^{1/x} W_{\alpha, \mu-i\tau}\left({1\over x}\right) W_{\beta, \nu+i\tau}\left({1\over x}\right) = 2  \int_{-\infty}^\infty  e^{-2i\tau y+ (\mu-\nu)y}  $$

$$\times  \int_0^{\pi/2}  \int_0^\infty e^{-x r^2 } r^{1-2(\alpha+\beta)}  \left({\cos\theta + e^y \sin\theta \over  e^y \cos\theta +\sin\theta }\right)^{\mu+\nu} $$

$$\times  K_{2(\mu+\nu)} \left( r \sqrt{ 1 + 2\cos\theta \sin\theta  \cosh y}\right)  {dr d\theta  dy\over \cos^{2\alpha}\theta\  \sin^{2\beta}\theta}.\eqno(1.13)$$
The inner integral with respect to $r$ is calculated as in (1.11). Hence, changing $1/x$ on $x$, we  derive from (1.13) the following product formula

$${ \Gamma\left(1/ 2-\alpha+\mu-i\tau\right)  \Gamma\left(1/ 2-\alpha-\mu+i\tau\right)   \Gamma\left(1/ 2-\beta+\nu+i\tau\right)  \Gamma\left(1/ 2-\beta-\nu-i\tau\right) \over \Gamma\left(1-\alpha-\beta + \mu+\nu\right) \Gamma\left(1-\alpha-\beta - \mu-\nu\right)}$$

$$\times  e^{x}\  W_{\alpha, \mu-i\tau}\left( x \right) W_{\beta, \nu+i\tau}\left( x\right) = \sqrt{x}  \int_{-\infty}^\infty  e^{-2i\tau y+ (\mu-\nu)y}  $$

$$\times  \int_0^{\pi/2}  {e^{x (1 + 2\cos\theta \sin\theta  \cosh y)/ 8}\over \sqrt{ 1 + 2\cos\theta \sin\theta  \cosh y} }\  \left({\cos\theta + e^y \sin\theta \over  e^y \cos\theta +\sin\theta }\right)^{\mu+\nu} $$

$$\times   W_{\alpha+\beta-1/2,  \mu+\nu}\left( {x\over 4} \ \left(1 + 2\cos\theta \sin\theta  \cosh y\right) \right)   { d\theta  dy\over \cos^{2\alpha}\theta\  \sin^{2\beta}\theta}.\eqno(1.14)$$
After simple changes of variables it can be written in a more symmetric form

$${ \Gamma\left(1/ 2-\alpha+\mu-i\tau\right)  \Gamma\left(1/ 2-\alpha-\mu+i\tau\right)   \Gamma\left(1/ 2-\beta+\nu+i\tau\right)  \Gamma\left(1/ 2-\beta-\nu-i\tau\right) \over \Gamma\left(1-\alpha-\beta + \mu+\nu\right) \Gamma\left(1-\alpha-\beta - \mu-\nu\right)}$$

$$\times  e^{x}\  W_{\alpha, \mu-i\tau}\left( x \right) W_{\beta, \nu+i\tau}\left( x\right) = \sqrt{x}  \int_{-\infty}^\infty  \int_{0}^{\pi/2}  e^{-2i\tau y+ (\mu-\nu)y}  $$

$$\times   {\left(1+\tan^2\theta\right)^{\alpha+1/4} \left(1+\cot^2\theta\right)^{\beta+1/4} \over \sqrt{ \tan\theta+\cot\theta + e^y+ e^{-y} }}\  \left({1 + e^y \tan\theta \over  e^y  +\tan\theta }\right)^{\mu+\nu} $$

$$\times \exp\left( {x \left(\tan\theta+\cot\theta + e^y+ e^{-y} \right) \over 8 (\tan\theta+\cot\theta)} \  \right) $$

$$\times   W_{\alpha+\beta-1/2,  \mu+\nu}\left( {x \left(\tan\theta+\cot\theta + e^y+ e^{-y} \right) \over 4 (\tan\theta+\cot\theta)} \  \right)  d\theta  dy$$

$$= \  \sqrt{x}  \int_{-\infty}^\infty  \int_{-\infty}^{\infty}  e^{-2i\tau y} \  {(1+ e^{2u})^{\alpha+\beta-1/2}\over e^{(2\beta- 1/2) u+ (\nu-\mu)y}}  \left({1 + e^{y+u} \over  e^{y}+ e^u}\right)^{\mu+\nu}  $$

$$\times \exp\left( {x\over 2}  \left(1 +  {\cosh y\over \cosh u} \right) \right)   W_{\alpha+\beta-1/2,  \mu+\nu}\left( x  \left(1 +  {\cosh y\over \cosh u} \right) \  \right)  {du dy\over \sqrt{\cosh u+\cosh y}},$$
i.e., finally,

$${ \Gamma\left(1/ 2-\alpha+\mu-i\tau\right)  \Gamma\left(1/ 2-\alpha-\mu+i\tau\right)   \Gamma\left(1/ 2-\beta+\nu+i\tau\right)  \Gamma\left(1/ 2-\beta-\nu-i\tau\right) \over \Gamma\left(1-\alpha-\beta + \mu+\nu\right) \Gamma\left(1-\alpha-\beta - \mu-\nu\right)}$$

$$\times  e^{x}\  W_{\alpha, \mu-i\tau}\left( x \right) W_{\beta, \nu+i\tau}\left( x\right) $$

$$ =  \sqrt{x}  \int_{-\infty}^\infty  \int_{-\infty}^{\infty}  e^{-2i\tau y} \  {(1+ e^{2u})^{\alpha+\beta-1/2}\over e^{(2\beta- 1/2) u+ (\nu-\mu)y}}  \left({1 + e^{y+u} \over  e^{y}+ e^u}\right)^{\mu+\nu}  $$

$$\times \exp\left( {x\over 2}  \left(1 +  {\cosh y\over \cosh u} \right) \right)   W_{\alpha+\beta-1/2,  \mu+\nu}\left( x  \left(1 +  {\cosh y\over \cosh u} \right) \  \right) {du dy\over \sqrt{\cosh u+\cosh y}}\eqno(1.15)$$
under conditions $x >0, \ \tau  \in \mathbb{R},\  \alpha, \beta < 1/2,  1-2\alpha >  2| \mu|,\  1-2\beta >  2|\nu|$. In particular, letting $\alpha,\beta > 0, \alpha+\beta = 1/2,\ \mu=\nu=0 $ and using (1.3), we obtain  the following representation for the product of Whittaker's functions of the pure imaginary index 

$$ |\Gamma\left(1/ 2-\alpha+i\tau\right) \Gamma\left(1/ 2-\beta+i\tau\right)|^2 \   W_{\alpha, i\tau}\left( x \right) W_{\beta, i\tau}\left( x\right) =  4 \sqrt{\pi } \ x e^{-x} $$

$$\times  \int_{0}^\infty  \int_{0}^{\infty}  \cos(2\tau y)\cosh\left(\left(2\beta- {1\over 2}\right)u\right) \ \exp\left( {x\over 2}  \left(1 +  {\cosh y\over \cosh u} \right) \right) $$

$$\times   K_0\left( {x\over 2}  \left(1 +  {\cosh y\over \cosh u} \right) \  \right) {du dy\over \sqrt{\cosh u}},\quad  x, \alpha,\beta > 0, \  \alpha+\beta= {1\over 2}.\eqno(1.16)$$
Finally, for $\alpha=\beta=1/4$ a double integral for the square of the Whittaker function $ W_{1/4, i\tau}^2\left( x \right)$

$$    W^2_{1/4, i\tau}\left( x \right)  =  {4 \sqrt{\pi } \ x e^{-x} \over |\Gamma\left(1/ 4+i\tau\right)|^4} \int_{0}^\infty  \int_{0}^{\infty}  \cos(2\tau y)\ \exp\left( {x\over 2}  \left(1 +  {\cosh y\over \cosh u} \right) \right) $$

$$\times   K_0\left( {x\over 2}  \left(1 +  {\cosh y\over \cosh u} \right) \  \right) {du dy\over \sqrt{\cosh u}},\quad  x > 0.\eqno(1.17)$$

Another approach to establish product formulas is to use  the Laplace integral  for the Whittaker function [3, formula (3.11)]

$$W_{\alpha, \mu-i\tau}(x)  = {e^{- x/2}\  x^{\mu-i\tau+1/2}   \over \Gamma(1/2-\alpha+\mu-i\tau) } \int_0^\infty  e^{-xs} s^{\mu-i\tau-\alpha-1/2}  (1+s)^{\mu-i\tau+\alpha-1/2}  ds \eqno(1.18)$$
under the convergence conditions $x >0, \   \mu+1/2 > \alpha$.  Hence we see that the  product $W_{\alpha, \mu-i\tau}(x) W_{\beta, \nu+i\tau}(x)$ can be treated as the Laplace transform of 
the corresponding convolution. Indeed, we have 
$$W_{\alpha, \mu-i\tau}(x) W_{\beta, \nu+i\tau}(x) = {e^{- x}\  x^{\mu+\nu+1}   \over \Gamma(1/2-\alpha+\mu-i\tau) \Gamma(1/2-\beta+\nu+i\tau)}$$

$$\times \int_0^\infty e^{- xs} \int_0^s t^{\mu-i\tau-\alpha-1/2}  (1+t)^{\mu-i\tau+\alpha-1/2} (s-t)^{\nu+i\tau-\beta-1/2} (1+s-t)^{\nu+i\tau+\beta-1/2} dt ds.\eqno(1.19)$$
The inner integral in (1.19) we write in the form

$$\int_0^s t^{\mu-i\tau-\alpha-1/2}  (1+t)^{\mu-i\tau+\alpha-1/2} (s-t)^{\nu+i\tau-\beta-1/2} (1+s-t)^{\nu+i\tau+\beta-1/2} dt$$

$$= s^{\mu+\nu-\alpha-\beta}   \int_0^1 t^{\mu-i\tau -\alpha-1/2}  (1+st)^{\mu-i\tau +\alpha-1/2} (1-t)^{\nu+i\tau-\beta-1/2} (1+s(1-t))^{\nu+i\tau+\beta-1/2} dt$$

$$= s^{\mu+\nu-\alpha-\beta}  (1+s)^{\nu+\mu+\beta+\alpha-1}  \int_0^1 t^{\mu-i\tau-\alpha-1/2}  \left(1-{s (1-t) \over 1+s} \right)^{\mu-i\tau+\alpha-1/2} (1-t)^{\nu+i\tau-\beta-1/2} $$

$$\times \left(1 - {st\over 1+s}\right)^{\nu+i\tau+\beta-1/2} dt = s^{\mu+\nu-\alpha-\beta}  (1+s)^{\nu+\mu+\beta+\alpha-1} $$

$$\times  \sum_{n=0}^\infty {(1/2-\mu+i\tau-\alpha)_n\  s^n\over n! \ (1+s)^n}  \int_0^1 t^{\mu-i\tau-\alpha-1/2} (1-t)^{\nu+i\tau+n-\beta-1/2}  $$

$$\times  \left(1 - {st\over 1+s}\right)^{\nu+i\tau+\beta-1/2}  dt,$$
where $(a)_n$ is the Pochhammer symbol and the interchange of the order of integration and summation is due to the dominated convergence theorem.  The latter integral is the Euler formula for the Gauss hypergeometric function (see Entry 2.2.6.1 in [2, Vol. I]), and we derive, combining with (1.19) the representation

$$ W_{\alpha, \mu-i\tau}(x) W_{\beta, \nu+i\tau}(x) = {e^{- x}\  x^{\mu+\nu+1}   \over \Gamma(\mu+\nu-\alpha-\beta+1) }$$

$$\times \int_0^\infty e^{- xs}  s^{\mu+\nu-\alpha-\beta}  (1+s)^{\nu+\mu+\beta+\alpha-1} $$

$$\times  \sum_{n=0}^\infty {(1/2-\mu+i\tau-\alpha)_n\  (1/2+\nu+i\tau-\beta)_n \over n! \ (\mu+\nu-\alpha-\beta+1)_n} \left({s\over 1+s}\right)^n $$

$$\times  \ {}_2F_1 \left( {1\over 2} -\nu-i\tau-\beta,\  {1\over 2} +\mu-i\tau-\alpha;\  \mu+\nu-\alpha-\beta+1+n;\  {s\over 1+s}\right) ds.\eqno(1.20)$$
The latter series is calculated via relation 6.7.1.2 in [2, Vol. III], and we establish the product formula in terms of the Appell $F_3$-function [1], [2]

$$ W_{\alpha, \mu-i\tau}(x) W_{\beta, \nu+i\tau}(x) = {e^{- x}\  x^{\mu+\nu+1}   \over \Gamma(\mu+\nu-\alpha-\beta+1) }$$

$$\times \int_0^\infty e^{- xs}  s^{\mu+\nu-\alpha-\beta}  (1+s)^{\nu+\mu+\beta+\alpha-1} \  F_3\left({1\over 2} -\nu-i\tau-\beta, \ {1\over 2} -\mu+i\tau-\alpha, \right.$$

$$\left.  {1\over 2}+\nu+i\tau-\beta,  \ {1\over 2} +\mu-i\tau-\alpha;\  \mu+\nu-\alpha-\beta+1; \ {s\over 1+s}, {s\over 1+s}\right) ds\eqno(1.21)$$
under conditions $x >0,\ \tau \in \mathbb{R},\ \mu+1/2 > \alpha,\ \nu +1/2 > \beta.$  Meanwhile, using  the auto-transformation formula for the Gauss hypergeometric function (cf. Entry 7.3.1.4 in [2, Vol. III]), we have 

$${}_2F_1 \left( {1\over 2} -\nu-i\tau -\beta,\  {1\over 2} +\mu-i\tau-\alpha;\  \mu+\nu-\alpha-\beta+1+n;\  {s\over 1+s}\right) = (s+1)^{-2\nu-2i\tau-n}$$

$$\times {}_2F_1 \left( {1\over 2} +\nu+i\tau-\beta+n,\  {1\over 2} + 2\nu+\mu+i\tau -\alpha+n;\  \mu+\nu-\alpha-\beta+1+n;\  {s\over 1+s}\right).$$
Thus equality (1.20) becomes

$$ W_{\alpha, \mu-i\tau}(x) W_{\beta, \nu+i\tau}(x) = {e^{- x}\  x^{\mu+\nu+1}   \over \Gamma(\mu+\nu-\alpha-\beta+1) }$$

$$\times \int_0^\infty e^{- xs}  s^{\mu+\nu-\alpha-\beta}  (1+s)^{\mu-\nu-2i\tau+\beta+\alpha-1} $$

$$\times  \sum_{n=0}^\infty {(1/2-\mu+i\tau-\alpha)_n\  (1/2+\nu+i\tau-\beta)_n \over n! \ (\mu+\nu-\alpha-\beta+1)_n} \left({s\over (1+s)^2}\right)^n $$

$$\times  \   {}_2F_1 \left(  {1\over 2} + 2\nu+\mu+i\tau -\alpha+n, \ {1\over 2} +\nu+i\tau-\beta+n;\  \mu+\nu-\alpha-\beta+1+n;\  {s\over 1+s}\right) ds.\eqno(1.22)$$
The series in (1.22) is calculated in [2, Vol.III], Entry 6.7.1.11 in terms of the Appell $F_1$-function, and we derive the product formula

$$ W_{\alpha, \mu-i\tau}(x) W_{\beta, \nu+i\tau}(x) = {e^{- x}\  x^{\mu+\nu+1}   \over \Gamma(\mu+\nu-\alpha-\beta+1) }$$

$$\times \int_0^\infty e^{- xs}  s^{\mu+\nu-\alpha-\beta}  (1+s)^{\mu-\nu-2i\tau+\beta+\alpha-1}  \  F_1\left({1\over 2} +\nu+i\tau-\beta, \ {1\over 2} -\mu+i\tau-\alpha, \right.$$

$$\left. \  2(\nu+\mu) ;  \ \mu+\nu-\alpha-\beta+1; \ {s\over 1+s},\  {s^2\over (1+s)^2}\right) ds.\eqno(1.23)$$
It is worth  mentioning an interesting particular case $\nu= - \mu$,  since relation 6.7.1.10 in [2, Vol. III] suggests a more simple form of the previous product formula. Precisely, we obtain 

$$ W_{\alpha, -\nu-i\tau}(x) W_{\beta, \nu+i\tau}(x) = {e^{- x}\  x   \over \Gamma(1-\alpha-\beta) }$$

$$\times \int_0^\infty e^{- xs}  s^{-\nu-\alpha-\beta}  (1+s)^{-2\nu-2i\tau+\beta+\alpha-1} $$

$$\times  \  {}_2F_1\left({1\over 2} +\nu+i\tau-\beta, \ {1\over 2} +\nu+i\tau-\alpha ;  \ 1-\alpha-\beta; \ {s^2+2s\over (1+s)^2}\right) ds.\eqno(1.24)$$
Using the same as in (1.9) transformation formula for the Gauss hypergeometric function and the parity of the Whittaker function with respect to its  second index, we write equality (1.24), accordingly,

$$ W_{\alpha, \nu+i\tau}(x) W_{\beta, \nu+i\tau}(x) = {e^{- x}\  x  \over \Gamma(1-\alpha-\beta) }$$

$$\times \int_0^\infty e^{- xs}  s^{-\alpha-\beta}  (1+s)^{\beta-\alpha} $$

$$\times  \  {}_2F_1\left(\ {1\over 2} +\nu+i\tau-\alpha,\ {1\over 2} -\nu-i\tau-\alpha;  \ 1-\alpha-\beta; \ - s^2 -2s\right) ds\eqno(1.25)$$
which coincides with (1.10) when $\nu=0.$

The main goal of this paper is to investigate the following class of index transformations

$$(F^\nu_{\alpha,\beta} f)(x)=  {e^x\over x} \int_{-\infty}^\infty   W_{\alpha, \nu+i\tau}(x) W_{\beta,\nu+i\tau}(x)$$

$$\times \Gamma\left({1\over 2} +\nu+i\tau-\alpha\right) \Gamma\left({1\over 2} -\nu-i\tau-\alpha\right)  f(\tau) d\tau,\quad x >0,\eqno(1.26)$$
involving three real parameters $\alpha, \beta, \nu$.  We will study  mapping properties of the operator $F^\nu_{\alpha,\beta}$ and prove inversion formulas for transformations $(F^\nu_{\alpha,\beta} f)(x)$  in suitable functional spaces. The case $\nu=0,\ \alpha +\beta =0$ was studied in [8]. Perhaps new features for this case will follow from the current investigation. Moreover, due to the value (1.3) this class of transformations covers integral operators with the product of the modified Bessel functions.  Finally in this section we observe that the key ingredient for the theory of these transformations is their relationship with the Laplace and Olevskii transforms [7] as we can see from (1.10) because the respective  Gauss hypergeometric function under the integral sign is the Olevskii kernel.  One more particular case relates to the equality $\nu=0, |\alpha-\beta|= 1/2$ since Entry 7.3.1.40 in [1, Vol. III] yields the expression of the Gauss hypergeometric function in terms of the Legendre function $P_{-1/2+i\tau}(z)$ being the kernel of the Mehler-Fock transform [6].  However, for the general case (1.26) we will need to involve the Olevskii transform of the  complex index $\nu+i\tau.$

\section{Olevskii transform of a  complex index with nonzero real part} 

Let us consider the following transformation of an arbitrary function $f$

$$(G_\nu f) (x)= {2^{2a-c-1}\  x^{c-1/2}\over  \Gamma(c)} \int_{-\infty}^\infty  \Gamma(a+\nu+i\tau) \Gamma(a-\nu-i\tau) \ $$

$$\times {}_2F_1\left( a + \nu+i\tau,\ a-\nu-i\tau; c ; - x^2 \right) f(\tau) d\tau,\quad x > 0,\eqno(2.1)$$
where $a, c, \nu$ are real parameters. The case $\nu=0$ corresponds to the classical Olevskii (also called Fourier-Jacobi ) transform [7], [9].  The key ingredient is the integral representation of the Olevskii kernel (see [9, formula (1.7)], relation 2.16.21.1 in [2, Vol. II]

$${}_2F_1\left( a + \nu+i\tau,\ a-\nu-i\tau; c ; - x^2 \right) = {2^{1-2a+c} x^{1-c} \Gamma(c)\over \Gamma(a+\nu+i\tau) \Gamma(a-\nu-i\tau) }$$

$$\times \int_0^\infty y^{2a-c} J_{c-1} (xy) K_{2(\nu+i\tau)} (y) dy,\quad x >0\eqno(2.2)$$
under conditions $ c >0,\ a >  |\nu| $.  Considering $f \in L_2\left(\mathbb{R}\right)$, we substitute the right-hand side of (2.2) in (2.1) and interchange the order of integration to obtain the equality

$$(G_\nu f)(x)=  \sqrt x \int_0^\infty y^{2a-c} J_{c-1} (xy)   \int_{-\infty}^\infty  K_{2(\nu+i\tau)} (y)  f(\tau)  d\tau dy,\quad x > 0.\eqno(2.3)$$
This equality implies the composition property of the Olevskii transform in terms of the Kontorovich-Lebedev and Hankel transforms [6].  The interchange is allowed due to the estimate for the Macdonald function (see [7, formula (2.121))

$$\left| K_{2(\nu+i\tau)} (y) \right| \le e^{- 2\delta |\tau|}  K_{2\nu} \left( y\cos\delta\right),\quad \delta \in \left(0,\ {\pi\over 2}\right)\eqno(2.4)$$
which justifies the absolute convergence of the iterated integral. Moreover,  its asymptotic behavior at zero and infinity [6] guarantees the property 

$$y^{2a-c-1/2}  \int_{-\infty}^\infty   K_{2(\nu+i\tau)} (y)  f(\tau) d\tau \in L_2\left(\mathbb{R}_+\right)$$
when $2a -c > 2 |\nu|$.  This condition guarantees the previous one $a > |\nu|$ since $ c > 0$.  Making a restriction $c \ge 1/2$,  we show that  integral (2.1) converges absolutely for each $x >0$.  Indeed, it  can be verified via (2.2), (2.4) with $\delta \in (0, \pi/2)$, employing  the Cauchy-Schwarz,   generalized Minkowski inequalities and  the uniform bound for  the Bessel function $|J_{c-1} (x)| < C x^{-1/2},\ x  > 0,\ c \ge 1/2$, where $C >0$ is an absolute constant. Precisely, we deduce

$$\int_{-\infty}^\infty  \left| \Gamma(a+\nu+i\tau) \Gamma(a-\nu-i\tau) \  {}_2F_1\left( a + \nu+i\tau,\ a-\nu-i\tau; c ; - x^2 \right) f(\tau) \right| d\tau$$

$$\le  2^{1-2a+c} x^{1-c} \Gamma(c)\  ||f||_{ L_2\left(\mathbb{R}\right)}   \left( \int_{-\infty}^\infty \left|  \int_0^\infty y^{2a-c} J_{c-1} (xy) K_{2(\nu+i\tau)} (y) dy\right|^2 d\tau \right)^{1/2}$$ 

$$\le  2^{1-2a+c} x^{1-c} \Gamma(c)\  ||f||_{ L_2\left(\mathbb{R}\right)} \int_0^\infty y^{2a-c} \left| J_{c-1} (xy)\right|  K_{2\nu} \left( y\cos\delta\right) dy $$ 

$$\times  \left( \int_{-\infty}^\infty   e^{- 4\delta |\tau|} d\tau \right)^{1/2} \le 2^{1-2a+c} C x^{1/2-c} \  \Gamma(c)\    ||f||_{ L_2\left(\mathbb{R}\right)} $$

$$\times  \int_0^\infty y^{2a-c-1/2}   K_{2\nu} \left( y\cos\delta\right) dy \left( \int_{-\infty}^\infty   e^{- 4\delta |\tau|} d\tau \right)^{1/2} < \infty,\ 2a -c > 2 |\nu|.$$

On the other hand,  the Parseval equality for the Hankel transform in $L_2$ yields

$$\int_0^\infty |(G_\nu f)(x)|^2 dx = \int_0^\infty  y^{4a-2c-1} \left|  \int_{-\infty}^\infty  K_{2(\nu+i\tau)} (y)\  f(\tau) d\tau \right|^2 dy.\eqno(2.5)$$
Now, using the representation of the Macdonald function in terms of the Fourier integral [7]

$$K_{2(\nu+i\tau)} (y) = {1\over 4}  \int_{-\infty}^\infty e^{-y\cosh (u/2)} e^{(\nu+i\tau)u} du,\ y > 0,\eqno(2.6)$$
and the Parseval equality for Fourier transform, equality (2.5) becomes

$$\int_0^\infty |(G_\nu f)(x)|^2 dx = {1\over 16} \int_0^\infty  y^{4a-2c-1} \left|  \int_{-\infty}^\infty e^{\nu u -y\cosh (u/2)} \hat{f} (u) du \right|^2 dy,\eqno(2.7)$$
where 

$$\hat{f} (u)= \int_{-\infty}^\infty f(\tau) e^{i\tau u} d\tau.\eqno(2.8)$$
Hence we have from (2.7)

$$\int_0^\infty |(G_\nu f)(x)|^2 dx = {1\over 16} \int_0^\infty  y^{4a-2c-1}  \int_{-\infty}^\infty   \int_{-\infty}^\infty e^{\nu (u+v)  -y(\cosh (u/2)+\cosh(v/2))} \hat{f} (u)  \overline{\hat{f}(v)} du dv dy$$

$$= {\Gamma(4a-2c) \over 16}  \int_{-\infty}^\infty   \int_{-\infty}^\infty {e^{\nu (u+v)}  \hat{f} (u)  \overline{\hat{f}(v)} \over (\cosh (u/2)+\cosh(v/2))^{4a-2c}}  du dv$$
and the Cauchy-Schwarz inequality for double integrals implies 

$$\int_0^\infty |(G_\nu f)(x)|^2 dx\le  {\Gamma(4a-2c) \over 16} \left(  \int_{-\infty}^\infty   \int_{-\infty}^\infty {e^{2\nu v}\  |\hat{f} (u)|^2   \over (\cosh (u/2)+\cosh(v/2))^{4a-2c}} du dv\right)^{1/2}$$

$$\times \left(  \int_{-\infty}^\infty   \int_{-\infty}^\infty {e^{2\nu u}\  |\hat{f}(v)|^2   \over (\cosh (u/2)+\cosh(v/2))^{4a-2c}} du dv\right)^{1/2}$$

$$\le   2^{4(a-1)- 2c}\ \Gamma(4a-2c)  \int_{-\infty}^\infty e^{2\nu v- (2a-c) |v|} dv   \int_{-\infty}^\infty  |\hat{f} (u)|^2 du $$

$$=   {\pi\  4^{2a- c-1} \ (2a-c)\  \Gamma(4a-2c)  \over  (2a-c)^2- 4\nu^2} \   ||f||^2_{L_2\left(\mathbb{R}\right)},\quad 2a-c > 2|\nu|.$$
Thus we derived the following $L_2$-norm estimate

$$  ||G_\nu||^2_{L_2\left(\mathbb{R}_+\right)} \le {\pi\  4^{2a- c-1}\  (2a-c)\ \Gamma(4a-2c) \over  (2a-c)^2- 4\nu^2} \   ||f||^2_{L_2\left(\mathbb{R}\right)},\quad 2a-c > 2|\nu|.\eqno(2.9)$$
Now, returning to (2.3), we invert the Hankel transform in $L_2$ to obtain

$$  \int_{-\infty}^\infty  K_{2(\nu+i\tau)} (y)  f(\tau)  d\tau =   y^{1-2a+c}  \int_0^\infty \sqrt{u} \  J_{c-1} (yu)\  (G_\nu f)(u) du,\quad y > 0,\eqno(2.10)$$
where the integral on the right-hand side converges in the mean square sense.  But it converges absolutely if $G_\nu \in L_2(\mathbb{R}_+)\cap L_1(\mathbb{R}_+)$ via  the uniform estimate above for the Bessel function when $c \ge 1/2$.   Hence multiplying both sides of (2.10) by $y^{\varepsilon -1} K_{2(\nu+ix)} (y), \ \varepsilon > 0,\ x \in \mathbb{R}$,  we integrate over $\mathbb{R}_+$ to obtain

$$ \int_0^\infty y^{\varepsilon -1} K_{2(\nu+ix)} (y) \int_{-\infty}^\infty  K_{2(\nu+i\tau)} (y)  f(\tau)  d\tau dy $$

$$ =  \int_0^\infty y^{\varepsilon +c-2a} K_{2(\nu+ix)} (y)  \int_0^\infty \sqrt{u} \  J_{c-1} (yu)\  (G_\nu f)(u) du dy.\eqno(2.11)$$
Then under  the condition $2a-c < 1/2-2|\nu|$ the interchange of the order of integration is possible on the latter iterated integral by Fubini's theorem.   Consequently, appealing again to Entry 2.16.21.1 in [2, Vol. II], we calculate the inner integral to derive the equality 

$$ \int_0^\infty y^{\varepsilon -1} K_{2(\nu+ix)} (y) \int_{-\infty}^\infty  K_{2(\nu+i\tau)} (y)  f(\tau)  d\tau dy $$
  
$$=  {2^{\varepsilon +c-2a-1} \over \Gamma(c)} \Gamma\left({\varepsilon\over 2} + c-a + \nu+ix\right) \Gamma\left({\varepsilon\over 2} + c-a - \nu-ix\right)$$

$$\times  \int_0^\infty u^{c-1/2}  \  {}_2F_1\left( {\varepsilon\over 2} + c-a + \nu+ix,\ {\varepsilon\over 2}+c-a-\nu-ix;\  c ; - u^2 \right) \ (G_\nu f)(u) du\eqno(2.12)$$
under the condition $c-a > |\nu|$.  Besides, it is possible to pass to the limit when $\varepsilon \to 0+$ under the integral sign on the right-hand side of (2.12). Indeed,  the integrand  is a continuous function of $\varepsilon$ on an interval $[ 0, \varepsilon_0],\ \varepsilon_0 > 0$ and the integral converges uniformly on this interval owing to the Weierstrass test.  Namely, it has when  $2a-c < 1/2-2|\nu|$

 $$\int_0^\infty y^{\varepsilon +c-2a} \left| K_{2(\nu+ix)} (y)\right|   \int_0^\infty \sqrt{u} \  \left|J_{c-1} (yu)\  (G_\nu f)(u) \right| du$$

$$\le \bigg[ A \int_0^1 y^{c-2a-1/2}  K_{2\nu} (y) dy + B \int_1^\infty  y^{\varepsilon_0+c-2a-1/2}  K_{2\nu} (y) dy \bigg]   \int_0^\infty | (G_\nu f)(u)|  du < \infty,$$
where $A, B$ are positive constants. Thus we obtain from (2.12) 

$$\lim_{\varepsilon \to 0+}  \int_0^\infty y^{\varepsilon -1} K_{2(\nu+ix)} (y) \int_{-\infty}^\infty  K_{2(\nu+i\tau)} (y)  f(\tau)  d\tau dy $$
  
$$=  {2^{c-2a-1} \over \Gamma(c)} \Gamma\left( c-a + \nu+ix\right) \Gamma\left( c-a - \nu-ix\right)$$

$$\times  \int_0^\infty u^{c-1/2}  \  {}_2F_1\left( c-a + \nu+ix,\ c-a-\nu-ix;\  c ; - u^2 \right) \ (G_\nu f)(u) du,\eqno(2.13)$$
where the integral on the right-hand side converges absolutely.   Taking into account other restrictions on parameters, we have  $c \ge  1/2,\ c-a > |\nu|, \ 2|\nu| < 2a-c < 1/2-2|\nu|$. 

Now, returning to (2.12), we will treat its left-hand side.  In fact, inequality (2.4) with  $\delta \in (0, \pi/2)$ and the condition  $f \in L_2\left(\mathbb{R}\right)$ allow us to interchange the order of integration for $|\nu| < \varepsilon/4$. Hence we appeal to Entry 2.16.33.2 in [2, Vol. II] to calculate the inner integral with respect to $y$, and we deduce

$$ \int_0^\infty y^{\varepsilon -1} K_{2(\nu+ix)} (y) \int_{-\infty}^\infty  K_{2(\nu+i\tau)} (y)  f(\tau)  d\tau dy $$

$$=  {2^{\varepsilon -3}\over \Gamma(\varepsilon)} \int_{-\infty}^\infty   \Gamma\left({\varepsilon\over 2} + 2\nu+i(x+\tau)\right) \Gamma\left({\varepsilon\over 2}  - 2\nu-i(x+\tau)\right) $$

$$\times  \Gamma\left({\varepsilon\over 2} +i(x-\tau)\right) \Gamma\left({\varepsilon\over 2}+ i(\tau-x)\right) f(\tau) d\tau,\ x \in \mathbb{R}.\eqno(2.14)$$
From now we are interested in the case $\nu\neq 0$. Hence we see that the right-hand side of (2.14) can be extended to  the interval $0 < \varepsilon < 4|\nu|$. Indeed, employing the familiar inequality for the gamma function $|\Gamma(z)| \le |\Gamma({\rm Re} z) |$ (see [6, Ch. 1]) and Entry 2.2.4.3 in [2, Vol. II], we have

$$ \int_{-\infty}^\infty  \left| \Gamma\left({\varepsilon\over 2} + 2\nu+i(x+\tau)\right) \Gamma\left({\varepsilon\over 2}  - 2\nu-i(x+\tau)\right) \right|$$

$$\times \left| \Gamma\left({\varepsilon\over 2} +i(x-\tau)\right) \Gamma\left({\varepsilon\over 2}+ i(\tau-x)\right) f(\tau) \right| d\tau $$

$$\le ||f||_{L_2(\mathbb{R})} \left| \Gamma\left({\varepsilon\over 2} + 2\nu\right) \Gamma\left({\varepsilon\over 2}  - 2\nu\right) \right| \left(\int_{-\infty}^\infty  \left| \Gamma\left({\varepsilon\over 2} +i(x-\tau)\right) \right|^4 d\tau \right)^{1/2} $$

$$=  \sqrt 2\ ||f||_{L_2(\mathbb{R})} \left| \Gamma\left({\varepsilon\over 2} + 2\nu\right) \Gamma\left({\varepsilon\over 2}  - 2\nu\right) \right| \left(\int_{0}^\infty  \left| \Gamma\left({\varepsilon\over 2} +i\tau)\right) \right|^4 d\tau \right)^{1/2} $$

$$=  ||f||_{L_2(\mathbb{R})} \left| \Gamma\left({\varepsilon\over 2} + 2\nu\right) \Gamma\left({\varepsilon\over 2}  - 2\nu\right) \right| {  2^{1-\varepsilon}\  \pi^{3/4} \ \Gamma^{3/2}(\varepsilon)\over \Gamma^{1/2} (1/2+ \varepsilon)}  < \infty .$$
Then, combining with (2.12), we derive the equality

$${1\over 2 \Gamma(\varepsilon)} \int_{-\infty}^\infty   \Gamma\left({\varepsilon\over 2} + 2\nu+i(x+\tau)\right) \Gamma\left({\varepsilon\over 2}  - 2\nu-i(x+\tau)\right) $$

$$\times  \Gamma\left({\varepsilon\over 2} +i(x-\tau)\right) \Gamma\left({\varepsilon\over 2}+ i(\tau-x)\right) f(\tau) d\tau $$

$$=  {2^{c-2a-\varepsilon+1} \over \Gamma(c)} \Gamma\left({\varepsilon\over 2} + c-a + \nu+ix\right) \Gamma\left({\varepsilon\over 2} + c-a - \nu-ix\right)  \int_0^\infty u^{c-1/2} $$

$$\times  {}_2F_1\left( {\varepsilon\over 2} + c-a + \nu+ix,\ {\varepsilon\over 2}+c-a-\nu-ix;\  c ; - u^2 \right) \ (G_\nu f)(u) du,\ x \in \mathbb{R}\eqno(2.15)$$
under  assumptions $c \ge 1/2,\  c-a > |\nu|, \ 2|\nu| < 2a-c < 1/2-2|\nu|,\  0 < \varepsilon < 4|\nu|$.   Denoting by $I_\nu(\varepsilon, x)$ the left-hand side of (2.15), we write it in the form (using the properties of the gamma function)

$$I_\nu(\varepsilon, x) = { \varepsilon \over 2 \Gamma(\varepsilon+1)} \int_{-\infty}^\infty   \Gamma\left({\varepsilon\over 2} + 2\nu+i(x+\tau)\right) \Gamma\left({\varepsilon\over 2}  - 2\nu-i(x+\tau)\right) $$

$$\times  {\Gamma\left(\varepsilon/ 2 +i(x-\tau)+1\right) \Gamma\left(\varepsilon/ 2+ i(\tau-x)+1\right)\over (\varepsilon/2)^2 + (\tau-x)^2}\   f(\tau) d\tau.\eqno(2.16)$$

{\bf Lemma 1}.  {\it Let $f \in L_2(\mathbb{R}),\  |\nu| < 1/2, \nu\neq 0$. Then for almost all $x \in \mathbb{R}$ the following relation holds}

$$\lim_{\varepsilon \to 0+} I_\nu(\varepsilon, x) = - {\pi^2\over 2}\  {f(x)\over (\nu+ix) \sin(2\pi(\nu+ix))}.\eqno(2.17)$$

\begin{proof}  The proof uses the fact (cf. [4]) that the limit relation

$$\lim_{\varepsilon \to 0+}  {1\over \pi }  \int_{-\infty}^\infty  {\varepsilon\ f(\tau) \over \varepsilon^2 + (\tau-x)^2}\  d\tau = f(x)$$
takes place for almost all $x \in \mathbb{R}$ if $f(\tau)/(1+\tau^2) \in L_1(\mathbb{R})$ which is true under the condition $f \in L_2(\mathbb{R})$ via the Cauchy-Schwarz inequality. So, writing (2.16) as follows

$$I_\nu(\varepsilon, x) =   {\varepsilon\over 2 \Gamma(\varepsilon+1)} \ \int_{-\infty}^\infty {f(\tau)\over (\varepsilon/2)^2 + (\tau-x)^2}$$

$$\times  \bigg[ \Gamma\left({\varepsilon\over  2} + 2\nu+i(x+\tau)\right) \Gamma\left({\varepsilon\over 2}  - 2\nu-i(x+\tau)\right)\bigg.$$

$$\times \Gamma\left({\varepsilon\over  2} +i(x-\tau)+1\right) \Gamma\left({\varepsilon\over 2}+ i(\tau-x)+1\right) $$

$$- \Gamma\left( 2(\nu+ix)\right) \Gamma\left(  - 2(\nu+ix)\right) \bigg]\ d\tau$$

$$+  {\varepsilon\  \Gamma\left( 2(\nu+ix)\right) \Gamma\left(  - 2(\nu+ix)\right) \over 2 \Gamma(\varepsilon+1)}  \int_{-\infty}^\infty {f(\tau) - f(x) \over (\varepsilon/2)^2 + (\tau-x)^2}\ d\tau $$ 

$$+  {\varepsilon\  f(x)\ \Gamma\left( 2(\nu+ix)\right) \Gamma\left(  - 2(\nu+ix)\right) \over 2 \Gamma(\varepsilon+1)}  \int_{-\infty}^\infty { d\tau \over (\varepsilon/2)^2 + (\tau-x)^2}$$ 
and taking into account the value of the integral

$$ {\varepsilon\over 2}  \int_{-\infty}^\infty { d\tau \over (\varepsilon/2)^2 + (\tau-x)^2} = \pi,\quad \varepsilon > 0,$$
we have

$$ \lim_{\varepsilon \to 0+} I_\nu(\varepsilon, x) =   \lim_{\varepsilon \to 0+}  {\varepsilon\over 2 \Gamma(\varepsilon+1)} \ \int_{-\infty}^\infty {f(\tau)\over (\varepsilon/2)^2 + (\tau-x)^2}$$

$$\times  \bigg[ \Gamma\left({\varepsilon\over  2} + 2\nu+i(x+\tau)\right) \Gamma\left({\varepsilon\over 2}  - 2\nu-i(x+\tau)\right)\bigg.$$

$$\times \Gamma\left({\varepsilon\over  2} +i(x-\tau)+1\right) \Gamma\left({\varepsilon\over 2}+ i(\tau-x)+1\right) $$

$$- \Gamma\left( 2(\nu+ix)\right) \Gamma\left(  - 2(\nu+ix)\right) \bigg]\ d\tau$$

$$+ \pi\   \Gamma\left( 2(\nu+ix)\right) \Gamma\left(  - 2(\nu+ix)\right) f(x).\eqno(2.18)$$
Thus it immediately implies (2.17) if we show that the limit on the right-hand side of (2.18) is zero.  To do this, we find

$$ {\varepsilon\over 2 \Gamma(\varepsilon+1)} \ \int_{-\infty}^\infty {|f(\tau)|\over (\varepsilon/2)^2 + (\tau-x)^2}$$

$$\times  \bigg|  \Gamma\left({\varepsilon\over  2} + 2\nu+i(x+\tau)\right) \Gamma\left({\varepsilon\over 2}  - 2\nu-i(x+\tau)\right)\bigg.$$

$$\times \Gamma\left({\varepsilon\over  2} +i(x-\tau)+1\right) \Gamma\left({\varepsilon\over 2}+ i(\tau-x)+1\right) $$

$$\bigg. - \Gamma\left( 2(\nu+ix)\right) \Gamma\left(  - 2(\nu+ix)\right)\bigg| \ d\tau$$

$$\le   {\varepsilon\ ||f||_{L_2(\mathbb{R})} \over 2 \Gamma(\varepsilon+1)} \ \left(\int_{-\infty}^\infty {1\over [(\varepsilon/2)^2 + (\tau-x)^2]^2}\right.$$

$$\times  \bigg|  \Gamma\left({\varepsilon\over  2} + 2\nu+i(x+\tau)\right) \Gamma\left({\varepsilon\over 2}  - 2\nu-i(x+\tau)\right)\bigg.$$

$$\times \Gamma\left({\varepsilon\over  2} +i(x-\tau)+1\right) \Gamma\left({\varepsilon\over 2}+ i(\tau-x)+1\right) $$

$$\left.\bigg. - \Gamma\left( 2(\nu+ix)\right) \Gamma\left(  - 2(\nu+ix)\right) \bigg|^2 \ d\tau\right)^{1/2}$$

$$=  { ||f||_{L_2(\mathbb{R})} \over \sqrt{\varepsilon/2}\  \Gamma(\varepsilon+1)} \ \left(\int_{-\infty}^\infty {dy \over ( 1+ y^2)^2}\right.$$

$$\times  \bigg| \Gamma\left({\varepsilon\over  2}\ \left( 1+ i y\right)  + 2(\nu+ix)\right) \Gamma\left({\varepsilon\over 2}\left( 1- i y\right)  - 2(\nu+ix)\right)\bigg.$$

$$\left.\times \left| \Gamma\left({\varepsilon\over  2}  \left( 1+ i y\right)+1\right) \right|^2   - \Gamma\left( 2(\nu+ix)\right) \Gamma\left(  - 2(\nu+ix)\right) \bigg|^2 \right)^{1/2}$$

$$\le   { ||f||_{L_2(\mathbb{R})} \over \sqrt{\varepsilon/2}\  \Gamma(\varepsilon+1)} \bigg[ \ \left(\int_{-\infty}^\infty {dy \over ( 1+ y^2)^2}\right.\bigg.$$

$$\times  \bigg| \Gamma\left({\varepsilon\over  2}\ \left( 1+ i y\right)  + 2(\nu+ix)\right)  \Gamma\left({\varepsilon\over 2}\left( 1- i y\right)  - 2(\nu+ix)\right)\bigg.$$

$$\left.\bigg. - \Gamma\left( 2(\nu+ix)\right) \Gamma\left(  - 2(\nu+ix)\right) \bigg|^2 \right)^{1/2}$$

$$+ \left(1+ \Gamma\left({\varepsilon\over 2} + 1\right) \right) \left(\int_{-\infty}^\infty {dy \over ( 1+ y^2)^2} \bigg| \Gamma\left({\varepsilon\over  2}\ \left( 1+ i y\right)  + 2(\nu+ix)\right)\bigg.\right.$$

$$\bigg.\left. \times  \Gamma\left({\varepsilon\over 2}\left( 1- i y\right)  - 2(\nu+ix)\right)\bigg|^2  \left| \Gamma\left({\varepsilon\over 2} \left( 1+ i y\right)+1\right)  -1 \right|^2  \right)^{1/2} \bigg].\eqno(2.19)$$
But, in turn, the use of the generalized Minkowski inequality yields 

$$\left(\int_{-\infty}^\infty {dy \over ( 1+ y^2)^2}\   \bigg| \Gamma\left({\varepsilon\over  2}\ \left( 1+ i y\right)  + 2(\nu+ix)\right)\bigg.\right.$$

$$\bigg.\left. \times  \Gamma\left({\varepsilon\over 2}\left( 1- i y\right)  - 2(\nu+ix)\right)\bigg|^2 \left| \Gamma\left({\varepsilon\over 2} \left( 1+ i y\right)+1\right)  -1 \right|^2  \right)^{1/2}$$

$$= \left(\int_{-\infty}^\infty {dy \over ( 1+ y^2)^2}\   \bigg| \Gamma\left({\varepsilon\over  2}\ \left( 1+ i y\right)  + 2(\nu+ix)\right)\bigg.\right.$$

$$\bigg.\left. \times  \Gamma\left({\varepsilon\over 2}\left( 1- i y\right)  - 2(\nu+ix)\right)\bigg|^2 \left| \int_0^\infty e^{-u} \left(u^{(1+iy)\varepsilon/2} -1 \right) du \right|^2  \right)^{1/2}$$

$$\le \int_0^\infty e^{-u}  \left(\int_{-\infty}^\infty {\left|u^{(1+iy)\varepsilon/2} -1 \right|^2 \over ( 1+ y^2)^2}\   \bigg| \Gamma\left({\varepsilon\over  2}\ \left( 1+ i y\right)  + 2(\nu+ix)\right)\bigg.\right.$$

$$\bigg.\left. \times  \Gamma\left({\varepsilon\over 2}\left( 1- i y\right)  - 2(\nu+ix)\right)\bigg|^2 dy  \right)^{1/2} du $$

$$=  {\varepsilon \over 2} \int_0^\infty e^{-u}  \left(\int_{-\infty}^\infty {\varepsilon  \over \varepsilon^2 + y^2}  \bigg| \Gamma\left({1\over  2}\ \left( \varepsilon + i y\right)  + 2(\nu+ix)\right)\bigg.\right.$$

$$\bigg.\left. \times  \Gamma\left({1\over 2}\left( \varepsilon - i y\right)  - 2(\nu+ix)\right)\bigg|^2\left|\int_1^u v^{(\varepsilon +iy)/2 -1} dv \right|^2 dy \right)^{1/2} du $$

$$\le  {\varepsilon  \over 2}  \left(\int_{-\infty}^\infty {\varepsilon  \over \varepsilon^2 + y^2}  \bigg| \Gamma\left({1\over  2}\ \left( \varepsilon + i y\right)  + 2(\nu+ix)\right)\bigg.\right.$$

$$\bigg.\left. \times  \Gamma\left({1\over 2}\left( \varepsilon - i y\right)  - 2(\nu+ix)\right)\bigg|^2 dy \right)^{1/2}  $$

$$\times \bigg[ \int_0^1 e^{-u}  \int_u^1 v^{\varepsilon/2 -1} dv  du + \int_1^\infty  e^{-u}  \int_1^u  v^{\varepsilon/2 -1} dv  du \bigg]$$

$$\le   {\varepsilon  \over 2}  \left(\int_{-\infty}^\infty {\varepsilon  \over \varepsilon^2 + y^2}  \bigg| \Gamma\left({1\over  2}\ \left( \varepsilon + i y\right)  + 2(\nu+ix)\right)\bigg.\right.$$

$$\bigg.\left. \times  \Gamma\left({1\over 2}\left( \varepsilon - i y\right)  - 2(\nu+ix)\right)\bigg|^2 dy \right)^{1/2}  \bigg[ \int_0^1 {1-e^{-v}\over v}\ dv  +   \int_1^\infty  e^{-v}  v^{2|\nu| -1} dv \bigg]\eqno(2.20)$$
since  $0 < \varepsilon < 4|\nu|$. Therefore we see from (2.19), (2.20)    that 

$$ { ||f||_{L_2(\mathbb{R})} \over \sqrt{\varepsilon/2}\  \Gamma(\varepsilon+1)}  \left(1+ \Gamma\left({\varepsilon\over 2} + 1\right) \right) \left(\int_{-\infty}^\infty {dy \over ( 1+ y^2)^2} \bigg| \Gamma\left({\varepsilon\over  2}\ \left( 1+ i y\right)  + 2(\nu+ix)\right)\bigg.\right.$$

$$\bigg. \left. \times  \Gamma\left({\varepsilon\over 2}\left( 1- i y\right)  - 2(\nu+ix)\right)\bigg|^2  \left| \Gamma\left({\varepsilon\over 2} \left( 1+ i y\right)+1\right)  -1 \right|^2  \right)^{1/2} $$

$$\le \sqrt{{\pi \varepsilon\over 2}}\  {||f||_{L_2(\mathbb{R})} \over \Gamma(\varepsilon+1)}  \left(1+ \Gamma\left({\varepsilon\over 2} + 1\right) \right) 
\bigg|  \Gamma\left({\varepsilon \over  2} + 2\nu\right) \Gamma\left({\varepsilon \over  2} - 2\nu\right) \bigg| $$

$$\times  \bigg[ \int_0^1 {1-e^{-v}\over v}\ dv  +   \int_1^\infty  e^{-v}  v^{2|\nu| -1} dv \bigg]= O(\sqrt{\varepsilon}),\ \varepsilon \to 0+.$$
Concerning the first term on the right-hand side of the inequality (2.19), we write by virtue of the Cauchy-Schwarz inequality and simple changes of variables

$$   { ||f||_{L_2(\mathbb{R})} \over \sqrt{\varepsilon/2}\  \Gamma(\varepsilon+1)}  \left(\int_{-\infty}^\infty {dy \over ( 1+ y^2)^2} \bigg| \Gamma\left({\varepsilon\over  2}\ \left( 1+ i y\right)  + 2(\nu+ix)\right)  \Gamma\left({\varepsilon\over 2}\left( 1- i y\right)  - 2(\nu+ix)\right)\bigg.\right.$$

$$\left.\bigg. - \Gamma\left( 2(\nu+ix)\right) \Gamma\left(  - 2(\nu+ix)\right) \bigg|^2 \right)^{1/2}$$

$$=  { ||f||_{L_2(\mathbb{R})} \over \sqrt{\varepsilon/2}\  \Gamma(\varepsilon+1)}  \left(\int_{-\infty}^\infty {dy \over ( 1+ y^2)^2} \right.$$

$$\times \bigg| {\Gamma\left(\varepsilon \left( 1+ i y\right)/2  + 2(\nu+ix)+1\right)  \Gamma\left(\varepsilon \left( 1- i y\right)/2  - 2(\nu+ix)+1\right)\over (\varepsilon \left( 1+ i y\right)/2  + 2(\nu+ix)) (\varepsilon \left( 1- i y\right)/2  - 2(\nu+ix))}\bigg.$$

$$\left.\bigg. - \Gamma\left( 2(\nu+ix)\right) \Gamma\left(  - 2(\nu+ix)\right) \bigg|^2 \right)^{1/2}$$

$$\le   {\varepsilon  ||f||_{L_2(\mathbb{R})} \over 2\ \Gamma(\varepsilon+1)}  \left(\int_{-\infty}^\infty { dy \over ((\varepsilon/2)^2 +y^2)^2\ [ (\varepsilon^2/4 -4\nu^2 + (y+2x)^2)^2 + 16\nu^2 (y+2x)^2 ]} \right.$$

$$\times \bigg| \Gamma\left(1+\varepsilon/2 +2\nu  + i(2x+y) \right)  \Gamma \left(1+\varepsilon/2 -2\nu  - i(2x+y) \right) \bigg.$$

$$\left.\bigg. - \Gamma\left(1+ 2(\nu+ix)\right) \Gamma\left( 1 - 2(\nu+ix)\right)\bigg|^2 \right)^{1/2}$$

$$+ \sqrt{{\varepsilon\over 2}}\  { ||f||_{L_2(\mathbb{R})}  |\Gamma\left( 2\nu\right) \Gamma\left(  - 2\nu\right) |\over  \Gamma(\varepsilon+1)}  \left(\int_{-\infty}^\infty {dy \over (\varepsilon^2/4 -4\nu^2 +y^2)^2 + 16\nu^2 y^2 }\right)^{1/2} $$

$$+ \sqrt{\varepsilon}\  {2 ||f||_{L_2(\mathbb{R})}\over \Gamma(\varepsilon+1)} \   |\Gamma\left( 2\nu\right) \Gamma\left(  - 2\nu\right)(\nu+ix) |$$

$$\times  \left(\int_{-\infty}^\infty { \varepsilon\  y^2 dy \over  ((\varepsilon/2)^2 +y^2)^2\ [ (\varepsilon^2/4 -4\nu^2 + (y+2x)^2)^2 + 16\nu^2 (y+2x)^2 ]} \right)^{1/2}. \eqno(2.21)$$
But since $ [ (\varepsilon^2/4 -4\nu^2 + (y+2x)^2)^2 + 16\nu^2 (y+2x)^2 ]^{-1}$ is bounded for all  $\nu\neq 0, \ y,x \in \mathbb{R}, \ 0 < \varepsilon < 4|\nu|$ and   the integral 

$$ \int_{-\infty}^\infty {dy \over (\varepsilon^2/4 -4\nu^2 +y^2)^2 + 16\nu^2 y^2 } = (4\nu^2 - \varepsilon^2/4)^{-1/2}$$

$$\times  \int_{-\infty}^\infty {dy \over (4\nu^2 - \varepsilon^2/4)\ (y^2-1)^2 + 16\nu^2 y^2 }\le  (4\nu^2 - \varepsilon^2/4)^{-3/2} \int_{|y| \le 1/2} {dy \over (1-y^2)^2}$$

$$+ {(4\nu^2 - \varepsilon^2/4)^{-1/2} \over 16\nu^2} \int_{1/2< |y| < \infty} {dy \over y^2} < \infty,$$
 we observe that the last two terms in (2.21) are $O(\sqrt\varepsilon),\ \varepsilon \to 0$. Finally,  the first term on the right-hand side of the inequality (2.21) can be  estimated as follows
 
$$  {\varepsilon  ||f||_{L_2(\mathbb{R})} \over 2\ \Gamma(\varepsilon+1)}  \left(\int_{-\infty}^\infty { dy \over ((\varepsilon/2)^2 +y^2)^2\ [ (\varepsilon^2/4 -4\nu^2 + (y+2x)^2)^2 + 16\nu^2 (y+2x)^2 ]} \right.$$

$$\times \bigg| \Gamma\left(1+\varepsilon/2 +2\nu  + i(2x+y) \right)  \Gamma \left(1+\varepsilon/2 -2\nu  - i(2x+y) \right) \bigg.$$

$$\left.\bigg. - \Gamma\left(1+ 2(\nu+ix)\right) \Gamma\left( 1 - 2(\nu+ix)\right)\bigg|^2 \right)^{1/2}$$

$$\le {M\over \sqrt{  \varepsilon} } \bigg[ \ \bigg| \Gamma\left(1+{\varepsilon\over 2} -2\nu \right)\bigg| \left( \int_{-\infty}^\infty { dy \over (1 +y^2)^2} \right.$$

$$\times \left.\bigg| \Gamma\left(1+ {\varepsilon\over 2}  \left(1 +iy\right)  +2(\nu + ix) \right)  - \Gamma \left(1+2(\nu+ix)\right) \bigg|^2 \right)^{1/2}$$

$$+ \bigg| \Gamma\left(1+2\nu \right)\bigg| \left(\int_{-\infty}^\infty { dy \over (1 +y^2)^2} \right.$$

$$\times \left.\bigg| \Gamma\left(1+{\varepsilon\over 2} \left(1 -iy\right)  -2(\nu+ ix) \right)  - \Gamma \left(1-2(\nu+ix)\right) \bigg|^2 \right)^{1/2}\bigg]$$

$$=  M \sqrt{{\varepsilon\over 2}}  \bigg[ \ \bigg| \Gamma\left(1+{\varepsilon\over 2} -2\nu \right)\bigg| \left( \int_{-\infty}^\infty { dy \over 1 +y^2}\bigg| \int_0^\infty e^{-u} u^{2(\nu + ix) } \left[ \int_1^u v^{\varepsilon \left(1 +iy\right)/2 -1} dv  \right] du \bigg|^2 \right)^{1/2}$$

$$+ \bigg| \Gamma\left(1+2\nu \right)\bigg| \left(\int_{-\infty}^\infty { dy \over 1 +y^2} \bigg| \int_0^\infty e^{-u} u^{- 2(\nu + ix) } \left[ \int_1^u v^{\varepsilon \left(1 -iy\right)/2-1} dv \right] du \bigg|^2 \right)^{1/2}\bigg],$$
where $M >0$ is an absolute constant.  Then, similarly to (2.20), the generalized Minkowski inequality yields

$$  M \sqrt{{\varepsilon\over 2}}  \bigg[ \ \bigg| \Gamma\left(1+{\varepsilon\over 2} -2\nu \right)\bigg| \left( \int_{-\infty}^\infty { dy \over 1 +y^2}\bigg| \int_0^\infty e^{-u} u^{2(\nu + ix) } \left[ \int_1^u v^{\varepsilon \left(1 +iy\right)/2 -1} dv  \right] du \bigg|^2 \right)^{1/2}$$

$$+ \bigg| \Gamma\left(1+2\nu \right)\bigg| \left(\int_{-\infty}^\infty { dy \over 1 +y^2} \bigg| \int_0^\infty e^{-u} u^{- 2(\nu + ix) } \left[ \int_1^u v^{\varepsilon \left(1 -iy\right)/2-1} dv \right] du \bigg|^2 \right)^{1/2}\bigg]$$

$$\le   M \sqrt{{\pi \varepsilon\over 2}}  \bigg[ \ \bigg| \Gamma\left(1+{\varepsilon\over 2} -2\nu \right)\bigg| \int_0^\infty e^{-u} u^{2\nu }  \bigg| \int_1^u v^{\varepsilon/2 -1} dv \bigg| du$$

$$+ \bigg| \Gamma\left(1+2\nu \right)\bigg|  \int_0^\infty e^{-u} u^{-2\nu }  \bigg|  \int_1^u v^{\varepsilon/2 -1} dv \bigg| du\bigg]$$

$$\le   M \sqrt{{\pi \varepsilon\over 2}}  \bigg[ \ \bigg| \Gamma\left(1+{\varepsilon\over 2} -2\nu \right)\bigg| \bigg[ \int_0^1 e^{-u} u^{2\nu } \log\left({1\over u}\right) du + \int_1^\infty  e^{-u} u^{1+2\nu + 4|\nu|}   du\bigg] $$

$$+ \bigg| \Gamma\left(1+2\nu \right)\bigg| \bigg[ \int_0^1 e^{-u} u^{-2\nu } \log\left({1\over u}\right) du + \int_1^\infty e^{-u} u^{1-2\nu+ 4|\nu| }  du\bigg] \bigg]$$

$$= O(\sqrt \varepsilon),\ \varepsilon \to 0+,\quad  |\nu| < {1\over 2}.$$
 Thus, returning to (2.18), we find that

$$ \lim_{\varepsilon \to 0+} I_\nu(\varepsilon, x) =  \pi\   \Gamma\left( 2(\nu+ix)\right) \Gamma\left(  - 2(\nu+ix)\right) f(x)$$
which coincides with (2.17) and completes the proof of Lemma 1.

\end{proof}

Returning to (2.13), (2.15) and appealing to Lemma 1, we are ready to summarize our results as the following theorem.

{\bf Theorem 1}.  {\it Let $f \in L_2(\mathbb{R}),\  \nu \in \mathbb{R} \backslash\{0\},\    c \ge 1/2,\  a < c-|\nu| $ such that $\ 2|\nu| < 2a-c < 1/2-2|\nu|$. Then the generalized Olevskii transform $(2.1)$ is a bounded operator $G_\nu: L_2(\mathbb{R}) \to L_2(\mathbb{R}_+)$, where its norm bound is given by $(2.9)$

$$||G_\nu|| \le  \sqrt \pi\  2^{2a- c-1}\  \left({(2a-c)\ \Gamma(4a-2c) \over  (2a-c)^2- 4\nu^2}\right)^{1/2} $$
and the integral $(2.1)$ converges absolutely for all $x >0$. Besides, if $G_\nu \in  L_2(\mathbb{R}_+) \cap  L_1(\mathbb{R}_+)$ then  for almost all $\tau \in \mathbb{R},  f(\tau)$ can be recovered by the inversion formula 

$$ f(\tau) = 2^{c-2a+1}\ {\Gamma\left( c-a + \nu+i\tau\right) \Gamma\left( c-a - \nu-i\tau\right)\over \pi\  \Gamma(c)  \Gamma\left( 2(\nu+i\tau)\right) \Gamma\left(  - 2(\nu+i\tau)\right) }$$

$$\times  \int_0^\infty x^{c-1/2}  \  {}_2F_1\left( c-a + \nu+i\tau,\ c-a-\nu-i\tau;\  c ; - x^2 \right)  (G_\nu f)(x) dx,\eqno(2.22)$$
where the integral converges absolutely.}

{\bf Remark 1}. The inversion formula for the classical Olevskii transform ( $\nu=0$) cannot be obtained directly from (2.22).  To recover it,  we need to adjust the proof of Lemma 1 for this case. 

Indeed, for particular examples of the integral transform (1.26) considered later,  it is worth  establishing the inversion formula for transformation (2.1) when $\nu=0.$  Recalling (2.16), letting $\nu=0$ and making simple substitutions, we write it in the form

$$I_0 (\varepsilon, x) = { \varepsilon \over 8 x \Gamma(\varepsilon+1)} \int_{-\infty}^\infty  \left| \Gamma\left({\varepsilon\over 2} +i(x+\tau) +1\right)\right|^2 $$

$$\times \left| \Gamma\left({\varepsilon\over 2} +i(x-\tau) +1\right)\right|^2   {f(\tau)+f(-\tau) \over (\varepsilon/2)^2 + (\tau-x)^2} \ {d\tau \over \tau},\quad x\neq 0.\eqno(2.23)$$

{\bf Lemma 2}.  {\it Let $f \in L_2( (-1,1); d\tau/\tau^2) \cap L_2( \mathbb{R}\backslash (-1,1); d\tau)$. Then for almost all $x \in \mathbb{R}$ the following relation holds}

$$\lim_{\varepsilon \to 0+} I_0 (\varepsilon, x) = \pi \left[ f(x)+ f(-x) \right] \left| \Gamma(2ix)\right|^2.\eqno(2.24)$$

\begin{proof}  Under the conditions of the lemma and for $x\neq0$, we have

$$  \lim_{\varepsilon \to 0+}   I_0(\varepsilon, x) =    \lim_{\varepsilon \to 0+}  {\varepsilon\over 8x \Gamma(\varepsilon+1)} \ \int_{-\infty}^\infty {f(\tau) +f(-\tau) \over  (\varepsilon/2)^2 + (\tau-x)^2}$$

$$\times  \bigg[ \left| \Gamma\left({\varepsilon\over 2} +i(x+\tau) +1\right)\right|^2  \left| \Gamma\left({\varepsilon\over 2} +i(x-\tau) +1\right)\right|^2  -  \left| \Gamma(2ix+1)\right|^2 \bigg] {d\tau \over \tau} $$

$$+  \pi \left[ f(x)+ f(-x) \right] \left| \Gamma(2ix)\right|^2.\eqno(2.25) $$ 
In the meantime,  the conditions of the lemma imply that $f(\tau)/\tau  \in L_2(\mathbb{R})$, we find via the Cauchy-Schwarz, Minkowski and generalized Minkowski inequalities 

$$  {\varepsilon\over 2} \ \int_{-\infty}^\infty { |f(\tau) +f(-\tau)| \over  (\varepsilon/2)^2 + (\tau-x)^2}$$

$$\times \left| \bigg[ \left| \Gamma\left({\varepsilon\over 2} +i(x+\tau) +1\right)\right|^2  \left| \Gamma\left({\varepsilon\over 2} +i(x-\tau) +1\right)\right|^2  -  \left| \Gamma(2ix+1)\right|^2 \bigg] \right|{d\tau \over |\tau|} $$

$$\le  {\varepsilon\over 2} \left( \int_{-\infty}^\infty  |f(\tau) +f(-\tau)|^2  {d\tau\over \tau^2} \right)^{1/2}  \bigg[ \left(\int_{-\infty}^\infty {d\tau  \over  [(\varepsilon/2)^2 + (\tau-x)^2]^2 }\right.$$

$$\times \left.  \bigg| \left| \Gamma\left({\varepsilon\over 2} +i(x+\tau) +1\right)\right|^2   -  \left| \Gamma(2ix+1)\right|^2 \bigg|^2 \right)^{1/2} $$

$$+ \left(\Gamma\left(1+{\varepsilon\over 2} \right) +1\right) \left(\int_{-\infty}^\infty {d\tau  \over  [(\varepsilon/2)^2 + (\tau-x)^2]^2 }\left| \Gamma\left({\varepsilon\over 2} +i(x+\tau) +1\right)\right|^4\right.$$

$$\times \left. \left|  \Gamma\left({\varepsilon\over 2} +i(x-\tau) +1\right)  -  1\right|^2 \right)^{1/2} \bigg] $$

$$\le  {1\over \sqrt{{\varepsilon/2}} }\left( \int_{-\infty}^\infty  |f(\tau) +f(-\tau)|^2  {d\tau\over \tau^2} \right)^{1/2}  \bigg[ \left(\int_{-\infty}^\infty {dy  \over  (1 + y^2)^2 }\right.$$

$$\times \left.  \bigg| \left| \Gamma\left({\varepsilon\over 2} (1+iy) +2ix+1\right)\right|^2   -  \left| \Gamma(2ix+1)\right|^2 \bigg|^2 \right)^{1/2} $$

$$+ \Gamma^2\left(1+{\varepsilon\over 2} \right) \left(\Gamma\left(1+{\varepsilon\over 2} \right) +1\right) \left(\int_{-\infty}^\infty {dy  \over  (1 +  y^2)^2 } \left|  \int_0^\infty e^{-u} \left(u^{(1-iy)\varepsilon/2} -1 \right) du\right|^2 \right)^{1/2} \bigg]$$

$$\le  {1\over \sqrt{{\varepsilon/2}} }\left( \int_{-\infty}^\infty  |f(\tau) +f(-\tau)|^2  {d\tau\over \tau^2} \right)^{1/2}  \bigg[ \left(\int_{-\infty}^\infty {dy  \over  (1 + y^2)^2 }\right.$$

$$\times \left.  \bigg| \left| \Gamma\left({\varepsilon\over 2} (1+iy) +2ix+1\right)\right|^2   -  \left| \Gamma(2ix+1)\right|^2 \bigg|^2 \right)^{1/2} $$

$$+ \Gamma^2\left(1+{\varepsilon\over 2} \right) \left(\Gamma\left(1+{\varepsilon\over 2} \right) +1\right) \int_0^\infty e^{-u}  \left(\int_{-\infty}^\infty { | u^{(1-iy)\varepsilon/2} -1 |^2   \over  (1 +  y^2)^2 } \ dy\right)^{1/2} du \bigg].\eqno(2.26)$$
But, in turn, for $0 < \varepsilon < 1$

$$\int_0^\infty e^{-u}  \left(\int_{-\infty}^\infty { | u^{(1-iy)\varepsilon/2} -1 |^2   \over  (1 +  y^2)^2 } \ dy\right)^{1/2} du \le  {\varepsilon\ \sqrt\pi \over 2} \int_0^\infty e^{-u}   \left| \int_1^u v^{\varepsilon/2-1} dv\right| du $$

$$\le  {\varepsilon\ \sqrt\pi \over 2} \bigg[ \int_0^1 {1- e^{-v}\over v}\ dv +    \int_1^\infty  e^{-v} dv \bigg] = O(\varepsilon),\ \varepsilon \to 0+.$$
Meanwhile, analogously,

$$\left(\int_{-\infty}^\infty {dy  \over  (1 + y^2)^2 }  \bigg| \left| \Gamma\left({\varepsilon\over 2} (1+iy) +2ix+1\right)\right|^2   -  \left| \Gamma(2ix+1)\right|^2 \bigg|^2 \right)^{1/2} $$

$$\le    \left(\Gamma\left(1+{\varepsilon\over 2} \right) +1\right)  \left(\int_{-\infty}^\infty {dy  \over  (1 + y^2)^2 }  \left| \int_0^\infty e^{-u} u^{2ix} \left(u^{\varepsilon (1+iy)/2}-1\right) du\right|^2 \right)^{1/2} $$

$$\le  {\varepsilon\ \sqrt\pi \over 2}   \left(\Gamma\left(1+{\varepsilon\over 2} \right) +1\right) \int_0^\infty e^{-u}   \left| \int_1^u v^{\varepsilon/2-1} dv\right| du =O(\varepsilon),\ \varepsilon \to 0+.$$
Consequently, combining with (2.25), (2.26), we get (2.24), completing the proof of Lemma 2.

\end{proof}

Thus, taking into account Lemma 2 and Theorem 1, the classical Olevskii transform (2.1)  $G_0 f$ under the corresponding conditions has the following inversion formula

$$ {f(\tau) + f(-\tau)\over 2} = 2^{c-2a}\ {|\Gamma\left( c-a +i\tau\right)|^2 \over \pi\  \Gamma(c)  |\Gamma\left( 2i\tau \right)|^2}$$

$$\times  \int_0^\infty x^{c-1/2}  \  {}_2F_1\left( c-a +i\tau,\ c-a-i\tau;\  c ; - x^2 \right)  (G_0 f)(x) dx.\eqno(2.27)$$
We observe that when $f$ is even,  it coincides with the respective inversion in [9] after changes of variables and functions.

\section{Mapping properties of the index transform (1.26)}

To proceed with the investigation of the transformation (1.26) we first estimate its kernel, using integral representation (1.25) under conditions $ x >0,\ \tau \in \mathbb{R}, \alpha+\beta \le 1/2$.   Indeed, recalling (2.2), the Gauss hypergeometric function is represented by the integral

$${}_2F_1\left(\ {1\over 2} +\nu+i\tau-\alpha,\ {1\over 2} -\nu-i\tau-\alpha;  \ 1-\alpha-\beta; \ - s^2 -2s\right)$$

$$ = {2^{1+\alpha-\beta}\  \Gamma( 1-\alpha-\beta ) \over \Gamma (1/2 +\nu+i\tau-\alpha) \Gamma (1/2-\nu-i\tau-\alpha)}$$

$$\times\  [s(s+2)]^{(\alpha+\beta)/2} \int_0^\infty y^{\beta-\alpha} J_{-\alpha-\beta}\left( y\ [s(s+2)]^{1/2}\right) K_{2(\nu+i\tau)}(y) dy,\eqno(3.1)$$
where $\alpha+\beta  \le 1/2,\  \alpha < 1/2-|\nu|$. Hence, employing the inequality (2.4) for the Macdonald function and  uniform bounds for the  Bessel function (see above), we derive

$$\left| {}_2F_1\left(\ {1\over 2} +\nu+i\tau-\alpha,\ {1\over 2} -\nu-i\tau-\alpha;  \ 1-\alpha-\beta; \ - s^2 -2s\right)\right| $$

$$ \le {2^{1+\alpha-\beta}\  \Gamma( 1-\alpha-\beta ) \over |\Gamma (1/2 +\nu+i\tau-\alpha) \Gamma (1/2-\nu-i\tau-\alpha)|}$$

$$\times  [s(s+2)]^{(\alpha+\beta)/2} e^{-2\delta |\tau|} \int_0^\infty y^{\beta-\alpha} \left|J_{-\alpha-\beta}\left( y\ [s(s+2)]^{1/2}\right) \right| K_{2\nu}(y\cos\delta) dy $$

$$ \le {2^{1+\alpha-\beta}\   C \Gamma( 1-\alpha-\beta )  e^{-2\delta |\tau|}   [s(s+2)]^{(\alpha+\beta)/2 -1/4} \over |\Gamma (1/2 +\nu+i\tau-\alpha) \Gamma (1/2-\nu-i\tau-\alpha)|}  $$
$$\bigg. \times \int_0^\infty y^{\beta-\alpha-1/2}  K_{2\nu}(y\cos\delta) dy = {e^{-2\delta |\tau|}  O \left( [s(s+2)]^{(\alpha+\beta)/2 -1/4}\right) \bigg] \over |\Gamma (1/2 +\nu+i\tau-\alpha) \Gamma (1/2-\nu-i\tau-\alpha)|} \eqno(3.2)$$
under conditions $|\nu| <  \min (1/2-\alpha,  (\beta +1/2-\alpha)/2 ),\  \alpha+\beta \le 1/2$ which yields as it is easily seen $|\nu| <  (\beta +1/2-\alpha)/2$.  Consequently,  recalling (1.25), the product of the Whittaker functions can be estimated accordingly, and we have  for $ x >0,\ \tau \in \mathbb{R}$

$$ \left| W_{\alpha,\nu+i\tau}(x) W_{\beta, \nu+i\tau}(x)\right| \le  { x e^{- x}\over \Gamma(1-\alpha-\beta) }$$

$$\times \int_0^\infty e^{- xs}  s^{-\alpha-\beta}  (1+s)^{\beta-\alpha} $$

$$\times  \ \left| {}_2F_1\left(\ {1\over 2} +\nu+i\tau-\alpha,\ {1\over 2} -\nu-i\tau-\alpha;  \ 1-\alpha-\beta; \ - s^2 -2s\right) \right|ds$$

$$\le  { x \  e^{-2\delta |\tau|-x}  \over  |\Gamma (1/2 +\nu+i\tau-\alpha) \Gamma (1/2-\nu-i\tau-\alpha)|}$$

$$\times \bigg[   O\left( \int_0^\infty e^{- xs}  s^{-(\alpha+\beta)/2- 1/4} \  (1+s)^{\beta-\alpha} (s+2)^{(\alpha+\beta)/2- 1/4}  ds\right)\bigg],\eqno(3.3) $$
where the integral on the right-hand side of the latter inequality can be written in terms of the   Tricomi function  $\Psi(a,b; z)$ [2, Vol. III]. Precisely, this gives

$$ \int_0^\infty e^{- xs}  s^{-(\alpha+\beta)/2- 1/4} \  (1+s)^{\beta-\alpha} (s+2)^{(\alpha+\beta)/2- 1/4}  ds $$

$$\le  \int_0^\infty e^{- xs}  s^{-(\alpha+\beta)/2- 1/4} \  (1+s)^{1/2- 2\alpha} ds$$

$$= \Gamma\left({3\over 4}- {\alpha+\beta\over 2}\right)  \Psi \left( {3\over 4} -  {\alpha+\beta \over 2},\  {9\over 4} - {\beta\over 2} -{5\alpha\over 2} ;\ x\right).$$
i.e. (3.3) becomes

$$ \left| W_{\alpha,\nu+i\tau}(x) W_{\beta, \nu+i\tau}(x)\right| \le { x \  e^{-2\delta |\tau|-x}  \over  |\Gamma (1/2 +\nu+i\tau-\alpha) \Gamma (1/2-\nu-i\tau-\alpha)|}$$

$$\times  O\left( \Psi \left( {3\over 4} -  {\alpha+\beta \over 2},\  {9\over 4} - {\beta\over 2} -{5\alpha\over 2} ;\ x\right)\right).\eqno(3.4)$$
The estimate (3.4) allows us to establish the composition structure of the index transform (1.26) in terms of the Laplace transform and generalized  Olevskii transform (2.1)  under suitable conditions on $f$ and to derive  its asymptotic bounds near zero and infinity.  In fact,  we have

{\bf Theorem 2}. {\it Let  $f \in L_1\left( \mathbb{R} ; e^{-2\delta |\tau|}  d\tau\right), \delta \in \left[0,  \pi/2\right)$ and parameters $\alpha, \beta, \nu \in \mathbb{R}$ such that $\alpha+\beta \le 1/2,\  |\nu| < (\beta -\alpha +1/2)/2$.  Then the transform $(1.26)$ exists as a continuous function on $\mathbb{R}_+$   and can be represented as a composition of the Laplace transform and generalized Olevskii transform, namely,
 
$$(F^\nu_{\alpha,\beta} f)(x)=   \int_0^\infty  t^{-\alpha-\beta}  e^{- x t}\  (G^\nu_{\alpha,\beta} f) (t) dt,\quad x >0,\eqno(3.5)$$
where

$$(G^\nu_{\alpha,\beta} f) (t)= {(1+t)^{\beta-\alpha}   \over \Gamma(1-\alpha-\beta) } \int_{-\infty}^\infty  \Gamma\left({1\over 2} +\nu+i\tau-\alpha\right) \Gamma\left({1\over 2} -\nu-i\tau-\alpha\right)$$

$$\times {}_2F_1\left(\ {1\over 2} +\nu+i\tau-\alpha,\ {1\over 2} -\nu-i\tau-\alpha;  \ 1-\alpha-\beta; \ - t^2 -2t\right) f(\tau) d\tau.\eqno(3.6)$$
Moreover, it has the following asymptotic bounds}

$$(F^\nu_{\alpha,\beta} f)(x) =  O\left(x^{(\alpha+\beta)/2-3/4}\right),\ x \to +\infty,\eqno(3.7)$$

$$(F^\nu_{\alpha,\beta} f)(x) =  O(1)+ O\left(x^{(5\alpha+\beta)/2- 5/4}\right),\ x \to 0,\   \beta + 5\alpha\neq {5\over 2},\eqno(3.8)$$

$$(F^\nu_{\alpha,\beta} f)(x) =  O\left(\log x\right),\ x \to 0,\  \beta + 5\alpha = {5\over 2}.\eqno(3.9)$$

 \begin{proof} Indeed, the proof is immediate from the dominated convergence and Fubini's theorems which guarantee the uniform convergence of the integral (1.26) on the interval $[x_0,\infty),\ x_0 >0$ and allow us to interchange the order of integration after substituting the integral representation (1.25).  Moreover,  by virtue of (3.4)

 $$ \left| (F^\nu_{\alpha,\beta} f)(x) \right| \le  O\left( \Psi \left( {3\over 4} -  {\alpha+\beta \over 2},\  {9\over 4} - {\beta\over 2} -{5\alpha\over 2} ;\ x\right)\right) ||f||_{L_1\left( \mathbb{R} ; e^{-2\delta |\tau|}  d\tau\right)}.$$
Since $x_0 > 0$ is arbitrary   we conclude that $(F^\nu_{\alpha,\beta} f)(x)$ is continuous on $\mathbb{R}_+$, having the composition representation (3.5). Finally, the asymptotic bounds (3.7)-(3.10) follow immediately from the asymptotic behavior of the Tricomi function near zero and infinity [2, Vol. III].  

 \end{proof}
 
 Further, Theorem 1 for the generalized   Olevskii transform (3.6) reads as follows :
 
  {\bf Theorem 3}.  {\it Let $f \in L_2(\mathbb{R}),\  \nu \in \mathbb{R} \backslash\{0\},\ \alpha, \beta \in \mathbb{R}$ such that  $\alpha+\beta \le 1/2$ and $|\nu| < (\beta-\alpha)/2 < 1/4-|\nu|,\ \beta < 1/2-|\nu|$.  Then the generalized Olevskii transform $(3.6)$ is a bounded operator 
 
 $$G^\nu_{\alpha,\beta} : L_2(\mathbb{R}) \to L_2\left(\mathbb{R}_+; \ \left[ t(t+2)\right]^{-\alpha-\beta} (t+1)^{1-2(\beta-\alpha)} dt\right)$$
and the inequality for norms takes place
$$  ||G^\nu_{\alpha,\beta} f ||^2_{L_2\left(\mathbb{R}_+;  \left[ t(t+2)\right]^{-\alpha-\beta} (t+1)^{1-2(\beta-\alpha)} dt\right)} \le {\pi\  (\beta-\alpha)\ \Gamma(2(\beta-\alpha)) \over (\beta-\alpha)^2- 4\nu^2} \   ||f||^2_{L_2\left(\mathbb{R}\right)},\eqno(3.11)$$
and the integral $(3.6)$ converges absolutely for all $x >0$. Besides, if 

$$G^\nu_{\alpha,\beta} f  \in  L_2\left(\mathbb{R}_+; \ \left[ t(t+2)\right]^{-\alpha-\beta} (t+1)^{1-2(\beta-\alpha)} dt\right)$$

$$ \cap\   L_1\left(\mathbb{R}_+; (t+1)^{1-\beta+\alpha} \left[ t (t+2)\right]^{-1/4-(\alpha+\beta)/2} dt\right),\eqno(3.12)$$ 
then  for almost all $\tau \in \mathbb{R},  f(\tau)$ can be recovered by the inversion formula 

$$ f(\tau) =  { \Gamma\left( 1/2-\beta+\nu+i\tau\right) \Gamma\left( 1/2-\beta-\nu-i\tau\right)\over \pi  \Gamma(1-\alpha-\beta)  \Gamma\left( 2(\nu+i\tau)\right) \Gamma\left(  - 2(\nu+i\tau)\right) }\  \int_0^\infty \left[ x (x+2)\right]^{-\alpha-\beta} (x+1)^{1+\alpha-\beta } $$

$$\times {}_2F_1\left( {1\over 2}+\nu+i\tau-\alpha,\  {1\over 2} -\nu -i\tau-\alpha;\  1-\nu-\alpha-\beta; - x^2 -2x\right) (G^\nu_{\alpha,\beta} f) (x) dx,\eqno(3.13)$$
where the integral converges absolutely.}

{\bf Remark 2}. When $\nu=0$ the corresponding inversion of (3.6) is guaranteed by additional assumptions of Lemma 2 and reads as follows (cf. (2.27))

$$ f(\tau) + f(-\tau) =  { |\Gamma\left( 1/2-\beta+i\tau\right) |^2\over \pi  \Gamma(1-\alpha-\beta) | \Gamma\left( 2i\tau\right) |^2 }\  \int_0^\infty \left[ x (x+2)\right]^{-\alpha-\beta} (x+1)^{1-(\beta-\alpha) } $$

$$\times {}_2F_1\left( {1\over 2}+i\tau-\alpha,\  {1\over 2} -i\tau-\alpha;\  1-\alpha-\beta; - x^2 -2x\right) (G^0_{\alpha,\beta} f) (x) dx.\eqno(3.14)$$

\section{Inversion formula}

The main result of this section is the inversion formula for the transformation (1.26).   To establish it, we return to (3.5) and  apply the Mellin transform [4]  over a vertical straight line in the complex plane $s$ to both its sides under the conditions of Theorems 2 and 3,  assuming $f \in L_2(\mathbb{R}) \subset  L_1\left( \mathbb{R} ; e^{-2\delta |\tau|}  d\tau\right),\ \delta \in \left(0,  \pi/2\right)$, and $ \alpha+\beta < 1/2,\  \beta < (\alpha+1/2)/6$. Then we  obtain via (3.5)

$$(F^\nu_{\alpha,\beta} f)^*(s) \equiv \int_0^\infty (F^\nu_{\alpha,\beta} f)(x)  x^{s-1} dx =  \Gamma(s)   \int_0^\infty  t^{-\alpha-\beta-s}  (G^\nu_{\alpha,\beta} f) (t) dt$$

$$ =   \Gamma(s)\  (G^\nu_{\alpha,\beta} f)^* (1-\alpha-\beta-s),\quad {\rm Re} (s) = \gamma \in \left(\beta-\alpha+{1\over 2},\   {3\over 4} -\ {\alpha+\beta\over 2} \right),\eqno(4.1)$$
 where the interchange of the order of integration is allowed due to the dominated convergence theorem.  Indeed,  from (3.6) via (3.2) we observe

 $$|(G^\nu_{\alpha,\beta} f) (t) | \le A_1\   (t+1)^{\beta-\alpha} [ t (t+2)]^{ (\alpha+\beta)/2- 1/4}\  ||f||_{L_2(\mathbb{R})} \eqno(4.2)$$
where $A_1 >0$ is a constant,  and the integral 

$$\int_0^\infty   t^{-(\alpha+\beta)/2-\gamma-1/4} (t+1)^{\beta-\alpha}  (t+2)^{ (\alpha+\beta)/2- 1/4}  dt < \infty$$
converges due to the interval for $\gamma$ in (4.1).    In the meantime, reciprocally, via the inverse Mellin transform in $L_2$ we have from (4.1)

 $$ (G^\nu_{\alpha,\beta} f) (t) = {1\over 2\pi i} \int_{1-\alpha-\beta-\gamma -i\infty}^{1-\alpha-\beta-\gamma+i\infty} {(F^\nu_{\alpha,\beta} f)^*(1-\alpha-\beta-s)\over \Gamma(1-\alpha-\beta-s)}\  t^{-s} ds.\eqno(4.3)$$
But the Parseval equality for the Mellin transform in $L_2$ [4] implies

$$\int_{0}^\infty |(G^\nu_{\alpha,\beta} f) (t)|^2 t^{-2(\alpha+\beta)-2\gamma +1} dt = {1\over 2\pi}  \int_{-\infty}^{\infty} \left|   {(F^\nu_{\alpha,\beta} f)^*(\gamma+iu) \over \Gamma(\gamma+ iu)} \right|^2 du.\eqno(4.4)$$
In fact, when  $ \gamma \in \left(\beta-\alpha+ 1/ 2,\   3/4 -\ (\alpha+\beta)/2 \right),\  \beta < (\alpha+1/2)/6$ the integral on the left-hand side of (4.4) converges.   Precisely, it has via (4.2)

$$ \int_{0}^\infty |(G^\nu_{\alpha,\beta} f) (t)|^2 t^{-2(\alpha+\beta)-2\gamma +1} dt = \left(\int_0^1 + \int_1^\infty \right) |(G^\nu_{\alpha,\beta} f) (t)|^2 t^{-2(\alpha+\beta)-2\gamma +1} dt  $$

$$\le B \int_0^1  t^{-\alpha-\beta-2\gamma +1/2} dt + C  \int_1^\infty  t^{2(\beta-\alpha-\gamma)} dt < \infty, $$
 where $B,C$ are positive constants.   Therefore, the convergence of the integral (4.3) is in the mean square sense with respect to the norm in  $L_2(\mathbb{R}_+;  t^{1-2(\alpha+\beta)-2\gamma } dt )$. Moreover,   the inequality (3.11) implies 
 
$$ G^\nu_{\alpha,\beta} f \in   L_2\left( \mathbb{R}_+; \left[ t(t+2)\right]^{-(\alpha+\beta)} (t+1)^{1-2(\beta-\alpha)}  d t\right) \cap\  L_2\left( \mathbb{R}_+;  t^{1-2(\alpha+\beta)-2\gamma  } dt \right),\eqno(4.5)$$
 where $\alpha+\beta < 1/2,\  \beta < (\alpha+1/2)/6,  \gamma \in \left(\beta-\alpha+1/ 2,\   3/4 - (\alpha+\beta)/2 \right)$.

 In the meantime, Entry 8.4.49.24 in [2, Vol. III]  gives the following Mellin-Barnes integral representation for the Gauss hypergeometric function

$$  {}_2F_1\left( {1\over 2}+\nu+i\tau-\alpha,\  {1\over 2} -\nu -i\tau-\alpha;\  1-\alpha-\beta; - x^2 -2x\right)  $$

$$\times \  {(x(x+2))^{- \alpha-\beta} \over \Gamma(1-\alpha-\beta) } = {1 \over 2\pi i }  \int_{\rho -i\infty}^{\rho +i\infty}  (x+1)^{2w} $$

$$\times  {\Gamma\left( 1/2 + \nu+\beta +i\tau+w \right) \Gamma\left( 1/ 2 -\nu+\beta-i\tau +w\right) \over \Gamma ( 1 + w)\ \Gamma( 1+\beta-\alpha +w)} dw,\eqno(4.7)$$
where $\rho > |\nu| -\beta-1/2$.  Meanwhile, it is possible to replace in (4.7)  the vertical line $(\rho-i\infty,\ \rho+i\infty)$ in the complex plane $w$ with  the left-hand loop $\cal{L}_{-\infty}$ with a  positive orientation and bounded imaginary part which encompasses   left-hand simple poles of gamma functions $w= -1/2-\beta \pm (\nu+i\tau) - n,\ n \in \mathbb{N}_0$. In fact, it is possible via analyticity property of the integrand and its exponential decay when $x > 0$ and ${\rm Re} (w) \to -\infty$ owing to the Stirling asymptotic formula for the gamma function.  Hence, inversion formula (3.13) can be written in the form

$$ f(\tau) =  { \Gamma\left( 1/2-\beta+\nu+i\tau\right) \Gamma\left( 1/2-\beta-\nu-i\tau\right)\over 2\pi^2 i \ \Gamma\left( 2(\nu+i\tau)\right) \Gamma\left(  - 2(\nu+i\tau)\right) }\  $$

$$\times  \int_{\cal{L}_{-\infty}} {\Gamma\left( 1/2 + \nu+\beta +i\tau+w \right) \Gamma\left( 1/ 2 -\nu+\beta-i\tau +w\right) \over \Gamma ( 1 + w)\ \Gamma( 1+\beta-\alpha + w)}$$

$$\times \int_0^\infty  (x+1)^{1+\alpha-\beta+2w }  (G^\nu_{\alpha,\beta} f) (x) dx dw,\eqno(4.8)$$
where the interchange of the order of integration is by virtue of the estimate (see (4.2))

$$ \int_{\cal{L}_{-\infty}} \bigg| {\Gamma\left( 1/2 + \nu+\beta +i\tau+w \right) \Gamma\left( 1/ 2 -\nu+\beta-i\tau +w\right) \over \Gamma ( 1 +w)\ \Gamma( 1+\beta-\alpha + w)}\bigg|$$

$$\times \int_0^\infty  (x+1)^{1+\alpha-\beta+2{\rm Re} (w) }  \left|(G^\nu_{\alpha,\beta} f) (x)\right| dx |dw| $$ 

$$\le  A_1 ||f||_{L_2(\mathbb{R})}  \int_{\cal{L}_{-\infty}} \bigg| {\Gamma\left( 1/2 + \nu+\beta +i\tau+w \right) \Gamma\left( 1/ 2 -\nu+\beta-i\tau +w\right) \over \Gamma ( 1 + w)\ \Gamma( 1+\beta-\alpha +w)}\bigg|$$

$$\times \int_0^\infty  (x+1)^{1+2{\rm Re} (w) }  [ x (x+2)]^{ (\alpha+\beta)/2- 1/4} dx |dw| $$ 

$$ \le  A_1 \Gamma\left( {3\over 4} + {\alpha+\beta\over 2} \right) ||f||_{L_2(\mathbb{R})}  \int_{\cal{L}_{-\infty}} \bigg| {\Gamma\left( 1/2 + \nu+\beta +i\tau+w \right) \Gamma\left( 1/ 2 -\nu+\beta-i\tau +w\right) \over \Gamma ( 1 + w)\ \Gamma( 1+\beta-\alpha +w) \Gamma(-3/4-(\alpha+\beta)/2-2w)}\bigg|$$

$$\times \left|\Gamma\left(- {3\over 2} - \alpha-\beta -2w \right)\right|  |dw|  < \infty$$ 
because the integrand in the latter integral  is $O( |w|^{ (\alpha+\beta)/2- 7/4}), \  {\rm Re} (w) \to -\infty$ and the  integral by $x$ converges when $\alpha+\beta > - 3/2, \ {\rm Re} (w) < - ( \alpha+\beta + 3/2)/2$.  From conditions of Theorem 3 we observe that $  |\nu| -\beta-3/2 <  - ( \alpha+\beta + 3/2)/2 <  |\nu| -\beta-1/2 $, i.e. $|\nu|- 3/4 < (\beta-\alpha)/2 < |\nu|+ 1/4$.  This means that the  left-hand loop $\hat{\cal{L}}_{-\infty}$  with  $ {\rm Re} (w) < - ( \alpha+\beta + 3/2)/2$ encompasses   all  simple poles of gamma functions $w= -1/2-\beta \pm (\nu+i\tau) - n,\ n \in \mathbb{N}$ except of two ones  $w= -1/2-\beta \pm (\nu+i\tau)$.  Thus, recalling (4.8), we write in the form

$$ f(\tau) =  { \Gamma\left( 1/2-\beta+\nu+i\tau\right) \Gamma\left( 1/2-\beta-\nu-i\tau\right)\over 2\pi^2 i \ \Gamma\left( 2(\nu+i\tau)\right) \Gamma\left(  - 2(\nu+i\tau)\right) }\  $$

$$\times  \int_{\hat{\cal{L}}_{-\infty}} {\Gamma\left( 1/2 + \nu+\beta +i\tau+w \right) \Gamma\left( 1/ 2 -\nu+\beta-i\tau +w\right) \over \Gamma ( 1 + w)\ \Gamma( 1+\beta-\alpha + w)}$$

$$\times \int_0^\infty  (x+1)^{1+\alpha-\beta+2w }  (G^\nu_{\alpha,\beta} f) (x) dx dw$$

$$+   { \Gamma\left( 1/2-\beta-\nu-i\tau\right)\over \pi \  \Gamma\left( 1/2-\alpha +\nu+i\tau\right) \Gamma\left(  - 2(\nu+i\tau)\right) }\  \int_0^\infty  (x+1)^{\alpha-3\beta+2(\nu+i\tau) }  (G^\nu_{\alpha,\beta} f) (x) dx$$

$$+ { \Gamma\left( 1/2-\beta+\nu+i\tau\right)\over \pi \  \Gamma\left( 1/2-\alpha -\nu-i\tau\right) \Gamma\left(   2(\nu+i\tau)\right) }\  \int_0^\infty  (x+1)^{\alpha-3\beta-2(\nu+i\tau) }  (G^\nu_{\alpha,\beta} f) (x) dx,\eqno(4.9)$$
where the convergence of the integrals will be clarified below.   Writing the power function  $ (x+1)^{1+\alpha-\beta+2w }$ in the form of the Laplace integral

$$ (x+1)^{1+\alpha-\beta+2w } = {1\over \Gamma(\beta-\alpha- 1- 2w )} \int_0^\infty e^{-(x+1)t} t^{\beta-\alpha- 2(1+w)}dt,$$
under the condition ${\rm Re} ( w) < (\beta-\alpha - 1)/2$,  we substitute it in the corresponding integral by $x$ in (4.9) and interchange the order of integration, having the expression

$$ \int_0^\infty (x+1)^{1+\alpha-\beta+2w }   (G^\nu_{\alpha,\beta} f) (x) dx =  {1\over \Gamma(\beta-\alpha- 1- 2w )} \int_0^\infty e^{-t}  t^{\beta-\alpha- 2(1+w)} $$

$$\times \int_0^\infty e^{-xt} (G^\nu_{\alpha,\beta} f) (x) dx  dt.$$
 This is, indeed, possible via (4.5) and  the Cauchy-Schwarz inequality, because
 
$$ \int_0^\infty e^{-t}  t^{\beta-\alpha- 2(1+ {\rm Re} (w))}  \int_0^\infty e^{-xt} \left| (G^\nu_{\alpha,\beta} f) (x)\right|  dx  dt $$

$$\le || G^\nu_{\alpha,\beta} f ||_{ L_2\left( \mathbb{R}_+;  x^{1-2(\alpha+\beta)-2\gamma  } dx \right)}  \int_0^\infty e^{-t}  t^{\beta-\alpha- 2(1+ {\rm Re} (w))}  \left(\int_0^\infty e^{-2 xt}  x^{2(\alpha+\beta+\gamma)-1 }  dx \right)^{1/2}  dt $$ 

$$= 2^{-\alpha-\beta-\gamma}\  \Gamma^{1/2} (2(\alpha+\beta+\gamma))\  \Gamma \left(-2( {\rm Re} (w)+\alpha)-\gamma-1\right) $$

$$\times  || G^\nu_{\alpha,\beta} f ||_{ L_2\left( \mathbb{R}_+;  x^{1-2(\alpha+\beta)-2\gamma  } dx \right)}  < \infty$$
under the condition $-\alpha-\beta < \gamma <  -2( {\rm Re} (w)+\alpha) -1$.  Consequently, recalling (4.3), (4.9) and the Mellin-Parseval equality, we find

$$  { \Gamma\left( 1/2-\beta+\nu+i\tau\right) \Gamma\left( 1/2-\beta-\nu-i\tau\right)\over 2\pi^2 i \ \Gamma\left( 2(\nu+i\tau)\right) \Gamma\left(  - 2(\nu+i\tau)\right) }\  $$

$$\times  \int_{\hat{\cal{L}}_{-\infty}} {\Gamma\left( 1/2 + \nu+\beta +i\tau+w \right) \Gamma\left( 1/ 2 -\nu+\beta-i\tau +w\right) \over \Gamma ( 1 + w)\ \Gamma( 1+\beta-\alpha + w)}$$

$$\times \int_0^\infty  (x+1)^{1+\alpha-\beta+2w }  (G^\nu_{\alpha,\beta} f) (x) dx dw$$

$$=  { \Gamma\left( 1/2-\beta+\nu+i\tau\right) \Gamma\left( 1/2-\beta-\nu-i\tau\right)\over 2\pi^2 i \ \Gamma\left( 2(\nu+i\tau)\right) \Gamma\left(  - 2(\nu+i\tau)\right) }\  $$

$$\times  \int_{\hat{\cal{L}}_{-\infty}} {\Gamma\left( 1/2 + \nu+\beta +i\tau+w \right) \Gamma\left( 1/ 2 -\nu+\beta-i\tau +w\right) \over \Gamma ( 1 + w)\ \Gamma( 1+\beta-\alpha + w) \Gamma(\beta-\alpha- 1- 2w )}$$

$$\times   \int_0^\infty e^{-t}  t^{\beta-\alpha- 2(1+w)}  \int_0^\infty e^{-xt} (G^\nu_{\alpha,\beta} f) (x) dx  dt dw $$

$$=  { \Gamma\left( 1/2-\beta+\nu+i\tau\right) \Gamma\left( 1/2-\beta-\nu-i\tau\right)\over \pi (2\pi i)^2  \ \Gamma\left( 2(\nu+i\tau)\right) \Gamma\left(  - 2(\nu+i\tau)\right) }\  $$

$$\times  \int_{\hat{\cal{L}}_{-\infty}} {\Gamma\left( 1/2 + \nu+\beta +i\tau+w \right) \Gamma\left( 1/ 2 -\nu+\beta-i\tau +w\right) \over \Gamma ( 1 + w)\ \Gamma( 1+\beta-\alpha + w) \Gamma(\beta-\alpha- 1- 2w )}$$

$$\times \int_{1-\alpha-\beta-\gamma -i\infty}^{1-\alpha-\beta-\gamma+i\infty} {(F^\nu_{\alpha,\beta} f)^*(1-\alpha-\beta-s) \over \Gamma(1-\alpha-\beta-s)}\     \Gamma (1-s)  \Gamma(\beta-\alpha+s - 2(1+w) ) ds dw$$

$$=  { \Gamma\left( 1/2-\beta+\nu+i\tau\right) \Gamma\left( 1/2-\beta-\nu-i\tau\right)\over 2\pi^2 i  \ \Gamma\left( 2(\nu+i\tau)\right) \Gamma\left(  - 2(\nu+i\tau)\right) } $$

$$\times \int_{1-\alpha-\beta-\gamma -i\infty}^{1-\alpha-\beta-\gamma+i\infty} {(F^\nu_{\alpha,\beta} f)^*(1-\alpha-\beta-s) \over \Gamma(1-\alpha-\beta-s)}\     \Phi^{\alpha,\beta}_{\nu+ i\tau} (s) \Gamma (1-s)  ds,\eqno(4.10)$$
where

$$\Phi^{\alpha,\beta}_{\nu+ i\tau} (s) =  {1\over 2\pi i} \int_{\hat{\cal{L}}_{-\infty}} {\Gamma\left( 1/2 + \beta+ \nu +i\tau+w \right) \Gamma\left( 1/ 2 +\beta-\nu -i\tau +w\right)  \over \Gamma ( 1 + w)\ \Gamma( 1+\beta-\alpha + w) \Gamma(\beta-\alpha- 1- 2w )}$$

$$\times \Gamma(\beta-\alpha+s - 2(1+w) ) dw,\eqno(4.11)$$
where the interchange of the order of integration is due to the absolute convergence of the iterated integral which can be shown, using the Cauchy-Schwarz inequality and since

$$ {\Gamma\left( 1/2 + \nu+\beta +i\tau+w \right) \Gamma\left( 1/ 2 -\nu+\beta-i\tau +w\right) \Gamma(\beta-\alpha+s - 2(1+w) ) \over \Gamma ( 1 + w)\ \Gamma( 1+\beta-\alpha + w) \Gamma(\beta-\alpha- 1- 2w )} $$

$$= O\left(  |w|^{\alpha+\beta+\gamma-2}\right),\quad \gamma < 1-\alpha-\beta,\  {\rm Re} (w) \to -\infty.$$
Furthermore, denoting by

$$\phi^{\alpha,\beta}_{\nu+ i\tau} (x) =  {1\over 2\pi i}  \int_{1-\alpha-\beta-\gamma -i\infty}^{1-\alpha-\beta-\gamma+i\infty} \Phi^{\alpha,\beta}_{\nu+ i\tau} (s)  x^{-s} ds,\quad x >0,\eqno(4.12)$$
we employ again the Mellin-Parseval equality and the fractional integration by parts to express the integral (4.10) as follows

$${1\over 2\pi i} \int_{1-\alpha-\beta-\gamma -i\infty}^{1-\alpha-\beta-\gamma+i\infty} {(F^\nu_{\alpha,\beta} f)^*(1-\alpha-\beta-s) \over \Gamma(1-\alpha-\beta-s)}\     \Phi^{\alpha,\beta}_{\nu+ i\tau} (s) \Gamma (1-s)  ds$$

$$= \int_0^\infty  \phi^{\alpha,\beta}_{\nu+ i\tau} (x)\ \left(I_-^{-\alpha-\beta} \left[ F^\nu_{\alpha,\beta} f\right] \right)(x) dx = \int_0^\infty  \left(I_{0+}^{-\alpha-\beta} \phi^{\alpha,\beta}_{\nu+ i\tau} \right) (t)\  \left(F^\nu_{\alpha,\beta} f\right)(t) dt,\eqno(4.13)$$
where $\left(I^\mu_{0+} g\right) (x),\  \left(I^\mu_- g\right) (x)$ are, correspondingly, the left-sided and right-sided Riemann-Liouville fractional integrals [6]

$$\left(I^\mu_{0+} g\right) (x) = {1\over \Gamma(\mu)}\int_0^x (x-t)^{\mu-1} g(t) dt,\quad x, \mu > 0,$$

$$\left(I^\mu_- g \right) (x) = {1\over \Gamma(\mu)}\int_x^\infty (t-x)^{\mu-1} g(t) dt,\quad x, \mu > 0.$$
Concerning other two integrals in (4.9), we have, accordingly, 

$$\int_0^\infty  (x+1)^{\alpha-3\beta\pm 2(\nu+i\tau) }  (G^\nu_{\alpha,\beta} f) (x) dx = {1\over \Gamma(3\beta-\alpha \mp 2(\nu+i\tau) )} \int_0^\infty e^{-t}  t^{3\beta-\alpha \mp 2(\nu+i\tau) -1} $$

$$\times \int_0^\infty e^{-xt} (G^\nu_{\alpha,\beta} f) (x) dx  dt,\eqno(4.14)$$
where $-1<  \alpha+\beta < 0,\ |\nu| < (3\beta-\alpha )/2$ and recalling  (4.5), we find the estimate

$$\int_0^\infty e^{-t}  t^{3\beta-\alpha \mp 2\nu -1}  \int_0^\infty e^{-xt} \left|(G^\nu_{\alpha,\beta} f) (x)\right| dx  dt $$

$$ \le \int_0^\infty e^{-t}  t^{3\beta-\alpha \mp 2\nu -1}  \left(\int_0^\infty e^{-2xt} \left[x(x+2)\right]^{\alpha+\beta} (x+1)^{2(\beta-\alpha)-1} dx\right)^{1/2} dx dt $$

$$\times \left|\left|G^\nu_{\alpha,\beta} f \right|\right|_{L_2\left( \mathbb{R}_+; \left[ x(x+2)\right]^{-(\alpha+\beta)} (x+1)^{1-2(\beta-\alpha)}  d x\right)}$$

$$ \le \Gamma^{1/2} (1+\alpha+\beta) \int_0^\infty e^{-t}  t^{3\beta-\alpha \mp 2\nu -1}  \Psi^{1/2}\left(1+\alpha+\beta, \ 1+4\beta; 2t\right) dt $$

$$\times \left|\left|G^\nu_{\alpha,\beta} f \right|\right|_{L_2\left( \mathbb{R}_+; \left[ x(x+2)\right]^{-(\alpha+\beta)} (x+1)^{1-2(\beta-\alpha)}  d x\right)},$$
where, in turn, the latter integral converges, using the asymptotic behavior of the Tricomi function, when $|\nu| < \min( (3\beta-\alpha )/2,\ (\beta-\alpha)/2)$.  Moreover, the composition of the Laplace transform and the generalized Olevskii transform (3.6) can be written in terms of the right-sided Riemann-Liouville fractional integral of the transformation (1.26). Precisely, we obtain, via (3.5),

$$\left( I_{-}^{-\alpha-\beta} \left[F^\nu_{\alpha,\beta} f\right] \right) (x) = {1\over \Gamma(-\alpha-\beta)}\int_x^\infty (y-x)^{-\alpha-\beta-1}  \left(F^\nu_{\alpha,\beta} f\right)(y) dy $$

$$=  {1\over \Gamma(-\alpha-\beta)}\int_x^\infty (y-x)^{-\alpha-\beta-1}  \int_0^\infty  t^{-\alpha-\beta}  e^{- y t}\  (G^\nu_{\alpha,\beta} f) (t) dt  dy $$

$$ = \int_0^\infty e^{-xt} (G^\nu_{\alpha,\beta} f) (t) dt.$$
Hence, returning to (4.14), we have owing to the fractional integration by parts 

$$\int_0^\infty  (x+1)^{\alpha-3\beta\pm 2(\nu+i\tau) }  (G^\nu_{\alpha,\beta} f) (x) dx =  {1\over \Gamma(3\beta-\alpha \mp 2(\nu+i\tau) )}$$

$$\times \int_0^\infty   \left(I_{0+}^{-\alpha-\beta} \left[e^{-x}  x^{3\beta-\alpha \mp 2(\nu+i\tau) -1} \right]\right) (t) \left(F^\nu_{\alpha,\beta} f\right) (t)dt.$$
Consequently, combining with (4.9)-(4.14), the inversion formula (3.13) takes the form

$$ f(\tau) =  \int_0^\infty {\cal K}_{\nu+i\tau}^{\alpha,\beta} (t) \left(F^\nu_{\alpha,\beta} f\right)(t) dt,\ \tau \in \mathbb{R},\eqno(4.15)$$
where 

$$ {\cal K}_{\nu+i\tau}^{\alpha,\beta} (x) =    { \Gamma\left( 1/2-\beta+\nu+i\tau\right) \Gamma\left( 1/2-\beta-\nu-i\tau\right)\over \pi  \ \Gamma\left( 2(\nu+i\tau)\right) \Gamma\left(  - 2(\nu+i\tau)\right) }  I_{0+}^{-\alpha-\beta} \left[  \phi^{\alpha,\beta}_{\nu+ i\tau} \right] (x)$$

$$+   { \Gamma\left( 1/2-\beta-\nu-i\tau\right)  \left( I_{0+}^{-\alpha-\beta} \left[  e^{-t} t^{3\beta-\alpha - 2(\nu+i\tau) -1} \right] \right) (x) \over \pi \  \Gamma\left( 1/2-\alpha +\nu+i\tau\right) \Gamma\left(  - 2(\nu+i\tau)\right) \Gamma(3\beta-\alpha - 2(\nu+i\tau) )}   $$

$$+ { \Gamma\left( 1/2-\beta+\nu+i\tau\right)  \left( I_{0+}^{-\alpha-\beta} \left[  e^{-t} t^{3\beta-\alpha + 2(\nu+i\tau) -1}\right] \right)(x)\over \pi \  \Gamma\left( 1/2-\alpha -\nu-i\tau\right) \Gamma\left(   2(\nu+i\tau)\right)  \Gamma(3\beta-\alpha + 2(\nu+i\tau) )},\ x >0.\eqno(4.16)$$
 The following inversion theorem summarizes the previous results of this section. 

{\bf Theorem 4}. {\it Let $ f\in L_2(\mathbb{R})$ and $\alpha, \beta, \nu \neq 0$ be real parameters, satisfying the following conditions

$$-1 < \alpha+\beta < 0,\quad \beta < {\alpha+ 1/2\over 6},$$

$$|\nu| < \min\left({3\beta -\alpha\over 2},\ {\beta-\alpha+1/2\over 2},\ {1\over 2} -\beta\right),$$

$$|\nu|  < {\beta-\alpha\over 2} < {1\over 4}- |\nu|.$$
Then for almost all $\tau \in \mathbb{R}$ the inversion formula for the transformation $(1.26)$ takes place  as the absolutely convergent integral $(4.15)$.  It involves  the kernel  $(4.16)$ ${\cal K}_{\nu+i\tau}^{\alpha,\beta} (x), \ x >0$, containing  the left-sided Riemann-Liouville fractional integrals of the power exponential functions and the the kernel $(4.12)$  $\phi^{\alpha,\beta}_{\nu+ i\tau} (x)$ which, in turn, is the inverse Mellin transform of  $\Phi^{\alpha,\beta}_{\nu+ i\tau} (s),\ {\rm Re}(s) = 1-\alpha-\beta-\gamma$ which is given by contour integral $(4.11)$  over the left-handed loop in the complex plane $w$ with restrictions}

$${\rm Re} (w) < \min \left( - {\alpha+\beta+3/2\over 2},\ {\beta-\alpha-1\over 2},\ - {2\alpha+1+\gamma\over 2}\right),$$

$$\max\left(-\alpha-\beta,\ \beta-\alpha+ {1\over 2}\right) < \gamma < \min\left( {3\over 4} - {\alpha+\beta\over 2},\ 1-\alpha-\beta\right).$$

\section{Representations for the kernel ${\cal K}_{\nu+i\tau}^{\alpha,\beta} (x)$. Particular cases}

Returning to (4.16),  we employ Entry 2.2.2.2 in [1] to calculate the fractional integral  $ \left( I_{0+}^{-\alpha-\beta} \left[  e^{-t} t^{3\beta-\alpha \pm 2(\nu+i\tau) -1}\right] \right) (x)$ in terms of the confluent hypergeometric function.  In fact, this gives

$$  \left(I_{0+}^{-\alpha-\beta} \left[  e^{-t} t^{3\beta-\alpha \pm 2(\nu+i\tau) -1}\right] \right) (x) =  {1\over \Gamma(-\alpha-\beta)} \int_0^x (x-t)^{-\alpha-\beta-1} e^{-t} t^{3\beta-\alpha \pm 2(\nu+i\tau) -1} dt$$

$$= x^{ 2(\beta-\alpha \pm (\nu+i\tau))-1} {\Gamma (3\beta-\alpha \pm 2(\nu+i\tau) )\over \Gamma( 2(\beta-\alpha \pm (\nu+i\tau)))}$$

$$\times  {}_1F_1 \left( 3\beta-\alpha \pm 2(\nu+i\tau) ;\  2(\beta-\alpha \pm (\nu+i\tau)) ;\ -x\right),\ x > 0.\eqno(5.1)$$
Meanwhile,   via (4.11), (4.12) we derive

$$ \left(I_{0+}^{-\alpha-\beta}  \phi^{\alpha,\beta}_{\nu+ i\tau}  \right) (x) =  {1\over 2\pi i}  \int_{1-\alpha-\beta-\gamma -i\infty}^{1-\alpha-\beta-\gamma+i\infty}  { \Phi^{\alpha,\beta}_{\nu+ i\tau} (s) \ \Gamma(1-s)\over \Gamma(1-\alpha-\beta-s)} x^{-s-\alpha-\beta} ds$$
and, accordingly, the Laplace transform equals to

$$\int_0^\infty e^{-xy} \left(I_{0+}^{-\alpha-\beta}  \phi^{\alpha,\beta}_{\nu+ i\tau} \right) (y) dy = {1\over 2\pi i}  \int_{1-\alpha-\beta-\gamma -i\infty}^{1-\alpha-\beta-\gamma+i\infty}   \Phi^{\alpha,\beta}_{\nu+ i\tau} (s)  \Gamma(1-s)\  x^{s+\alpha+\beta-1} ds$$

$$=   {x^{\alpha+\beta} \over 2\pi i} \int_{\hat{\cal{L}}_{-\infty}} {\Gamma\left( 1/2 + \beta+ \nu +i\tau+w \right) \Gamma\left( 1/ 2 +\beta-\nu -i\tau +w\right) \over \Gamma ( 1 + w)\ \Gamma( 1+\beta-\alpha + w) \Gamma(\beta-\alpha- 1- 2w )}$$

$$\times   {1\over 2\pi i}  \int_{\alpha+\beta+\gamma -i\infty}^{\alpha+\beta+\gamma+i\infty}   \Gamma(s)  \Gamma(\beta-\alpha-s - 2w-1 ) x^{-s} ds dw$$

$$=   {x^{\alpha+\beta} \over 2\pi i} \int_{\hat{\cal{L}}_{-\infty}} {\Gamma\left( 1/2 + \beta+ \nu +i\tau+w \right) \Gamma\left( 1/ 2 +\beta-\nu -i\tau +w\right) \over \Gamma ( 1 + w)\ \Gamma( 1+\beta-\alpha + w) }$$

$$\times  (x+1)^{2w+1+\alpha-\beta } dw =   \int_0^\infty e^{-xy}  \left(I_{0+}^{-\alpha-\beta} g\right)  (y) dy,\eqno(5.2)$$
where the function $g(x)$ is such that its Laplace transform is equal to

 $$ {1 \over 2\pi i}  \int_{\hat{\cal{L}}_{-\infty}} {\Gamma\left( 1/2 + \beta+ \nu +i\tau+w \right) \Gamma\left( 1/ 2 +\beta-\nu -i\tau +w\right) \over \Gamma ( 1 + w)\ \Gamma( 1+\beta-\alpha + w) }$$

$$\times  (x+1)^{2w+1+\alpha-\beta } dw,\eqno(5.3) $$
and  the interchange of the order of integration is allowed under conditions above.  It turns out that via injectivity of the composition  of the Laplace transform and fractional integral $g(x) = \phi^{\alpha,\beta}_{\nu+ i\tau} (x)$.  Therefore from (5.3) we have

$$\int_0^\infty e^{-xy} \phi^{\alpha,\beta}_{\nu+ i\tau} (y) dy =  {1 \over 2\pi i}  \int_{\hat{\cal{L}}_{-\infty}} {\Gamma\left( 1/2 + \beta+ \nu +i\tau+w \right) \Gamma\left( 1/ 2 +\beta-\nu -i\tau +w\right) \over \Gamma ( 1 + w)\ \Gamma( 1+\beta-\alpha + w) }$$

$$\times  (x+1)^{2w+1+\alpha-\beta } dw,\ x >0.\eqno(5.4) $$
The integral (5.3) can be expressed in terms of the Gauss hypergeometric function, returning to (4.7) and using the definition of the loop $\hat{\cal{L}}_{-\infty}$. Indeed, we obtain owing to the residue Cauchy theorem

$${1 \over 2\pi i}  \int_{\hat{\cal{L}}_{-\infty}} {\Gamma\left( 1/2 + \beta+ \nu +i\tau+w \right) \Gamma\left( 1/ 2 +\beta-\nu -i\tau +w\right) \over \Gamma ( 1 + w)\ \Gamma( 1+\beta-\alpha + w) }$$

$$\times  (x+1)^{2w+1+\alpha-\beta } dw $$

$$= (x+1)^{\alpha-3\beta} \bigg[{ (x+1)^{2(\nu+i\tau)}\  \Gamma(2(\nu+i\tau)) \over \Gamma(1/2-\beta+\nu+i\tau) \Gamma(1/2-\alpha+\nu+i\tau) } $$

$$\times    {}_2F_1\left( {1\over 2}-\nu-i\tau+\alpha,\  {1\over 2} -\nu -i\tau+\beta;\  1- 2(\nu+i\tau);  {1\over (x+1)^2}\right)\bigg.$$

$$+ { (x+1)^{- 2(\nu+i\tau)}\  \Gamma(-2(\nu+i\tau)) \over \Gamma(1/2-\beta-\nu-i\tau) \Gamma(1/2-\alpha-\nu-i\tau) } $$

$$\times    {}_2F_1\left( {1\over 2}+\nu+i\tau+\alpha,\  {1\over 2} +\nu +i\tau+\beta;\  1+ 2(\nu+i\tau);  {1\over (x+1)^2}\right)$$

$$-  { (x+1)^{ - 2(\nu+i\tau)}\  \Gamma(-2(\nu+i\tau)) \over \Gamma(1/2-\beta-\nu-i\tau) \Gamma(1/2-\alpha-\nu-i\tau) } -  { (x+1)^{2(\nu+i\tau)}\  \Gamma(2(\nu+i\tau)) \over \Gamma(1/2-\beta+\nu+i\tau) \Gamma(1/2-\alpha+\nu+i\tau) }\bigg].$$
Consequently, combining with (5.4) and appealing to Entries 2.2.2.1, 3.35.1.15 in [2, Vol. V], we invert the Laplace transform to obtain the value of the kernel $\phi^{\alpha,\beta}_{\nu+ i\tau} (x)$, namely,

$$ \phi^{\alpha,\beta}_{\nu+ i\tau} (x) =  e^{-x} x^{3\beta-\alpha -1} \bigg[-  { x^{ 2(\nu+i\tau)} \Gamma(-2(\nu+i\tau)) \over \Gamma(1/2-\beta-\nu-i\tau) \Gamma(1/2-\alpha-\nu-i\tau)  \Gamma (3\beta-\alpha+ 2(\nu+i\tau))}\bigg.$$

$$-  {x^{ -2(\nu+i\tau)}  \Gamma(2(\nu+i\tau)) \over \Gamma(1/2-\beta+\nu+i\tau) \Gamma(1/2-\alpha+\nu+i\tau) \Gamma (3\beta-\alpha- 2(\nu+i\tau))}$$

$$+ { x^{- 2(\nu+i\tau)}\  \Gamma(2(\nu+i\tau)) \over \Gamma(1/2-\beta+\nu+i\tau) \Gamma(1/2-\alpha+\nu+i\tau) \Gamma (3\beta-\alpha- 2(\nu+i\tau))} $$

$$\times    {}_2F_3\left( {1\over 2}-\nu-i\tau+\alpha,\  {1\over 2} -\nu -i\tau+\beta;\  1- 2(\nu+i\tau), \right.$$

$$\left.  {3\beta-\alpha\over 2}-\nu-i\tau,\  {1+3\beta-\alpha\over 2}-\nu-i\tau;  {x^2\over 4}\right)$$

$$+ { x^{ 2(\nu+i\tau)}\  \Gamma(-2(\nu+i\tau)) \over \Gamma(1/2-\beta-\nu-i\tau) \Gamma(1/2-\alpha-\nu-i\tau) \Gamma (3\beta-\alpha+ 2(\nu+i\tau))} $$

$$\times    {}_2F_3\left( {1\over 2}+\nu+i\tau+\alpha,\  {1\over 2} +\nu +i\tau+\beta;\  1+2(\nu+i\tau), \right.$$

$$\left.  {3\beta-\alpha\over 2}+\nu+i\tau,\  {1+3\beta-\alpha\over 2}+\nu+i\tau;  {x^2\over 4}\right)\bigg],\ x >0.\eqno(5.5)$$
Hence, recalling (5.4), we rewrite the kernel $ {\cal K}_{\nu+i\tau}^{\alpha,\beta} (x)$ in the form

 $$ {\cal K}_{\nu+i\tau}^{\alpha,\beta} (x) =  {  \Gamma(1/2-\beta-\nu-i\tau) \over \pi \Gamma(- 2(\nu+i\tau))  \Gamma(1/2-\alpha+\nu+i\tau) \Gamma (3\beta-\alpha- 2(\nu+i\tau))} $$

$$\times \left(I_{0+}^{-\alpha-\beta} \left[  e^{-t} t^{3\beta-\alpha - 2(\nu+i\tau) -1}   {}_2F_3\left( {1\over 2}-\nu-i\tau+\alpha,\  {1\over 2} -\nu -i\tau+\beta;\  1- 2(\nu+i\tau), \right.\right.\right.$$

$$\left. \left.\left. {3\beta-\alpha\over 2}-\nu-i\tau,\  {1+3\beta-\alpha\over 2}-\nu-i\tau;  {t^2\over 4}\right)\right] \right) (x) $$

$$+ {   \Gamma(1/2-\beta+\nu+i\tau) \over \pi \Gamma(2(\nu+i\tau))  \Gamma(1/2-\alpha-\nu-i\tau) \Gamma (3\beta-\alpha+ 2(\nu+i\tau))} $$

$$\times  \left(I_{0+}^{-\alpha-\beta} \left[  e^{-t} t^{3\beta-\alpha + 2(\nu+i\tau) -1}   {}_2F_3\left( {1\over 2}+\nu+i\tau+\alpha,\  {1\over 2} +\nu +i\tau+\beta;\  1+2(\nu+i\tau), \right.\right.\right.$$

$$\left. \left.\left. {3\beta-\alpha\over 2}+\nu+i\tau,\  {1+3\beta-\alpha\over 2}+\nu+i\tau;  {t^2\over 4}\right) \right] \right) (x),\ x >0.\eqno(5.6)$$
The fractional integrals in (5.6) we calculate directly in terms of the double hypergeometric series [10] via the definition of the hypergeometric ${}_2F_3$-function and using (5.1). In fact, we derive

$${  \Gamma(1/2-\beta-\nu-i\tau) \over \pi \Gamma(- 2(\nu+i\tau))  \Gamma(1/2-\alpha+\nu+i\tau) \Gamma (3\beta-\alpha- 2(\nu+i\tau))} $$

$$ \times \left(I_{0+}^{-\alpha-\beta} \left[  e^{-t} t^{3\beta-\alpha - 2(\nu+i\tau) -1}   {}_2F_3\left( {1\over 2}-\nu-i\tau+\alpha,\  {1\over 2} -\nu -i\tau+\beta;\  1- 2(\nu+i\tau), \right.\right.\right.$$

$$\left. \left.\left. {3\beta-\alpha\over 2}-\nu-i\tau,\  {1+3\beta-\alpha\over 2}-\nu-i\tau;  {t^2\over 4}\right)\right] \right) (x) $$

$$={  \Gamma(1/2-\beta-\nu-i\tau) \over \pi \Gamma(- 2(\nu+i\tau))  \Gamma(1/2-\alpha+\nu+i\tau) \Gamma( 2(\beta-\alpha - \nu-i\tau))} $$

$$\times \sum_{n=0}^\infty {(1/2-\nu-i\tau+\alpha)_n\  (1/2-\nu-i\tau+\beta)_n \  x^{ 2(\beta-\alpha - \nu-i\tau+n)-1} \over (1-2(\nu+i\tau))_n\  (2(\beta-\alpha - \nu-i\tau))_{2n}\   n!}$$

$$\times   {}_1F_1 \left( 3\beta-\alpha - 2(\nu+i\tau-n) ;\  2(\beta-\alpha - \nu-i\tau+n) ;\ -x\right)$$

$$={  x^{ 2(\beta-\alpha - \nu-i\tau)-1} \Gamma(1/2-\beta-\nu-i\tau) \over \pi \Gamma(- 2(\nu+i\tau))  \Gamma(1/2-\alpha+\nu+i\tau) \Gamma( 2(\beta-\alpha - \nu-i\tau))} $$

$$\times \sum_{n,k =0}^\infty  { (- 1)^k x^{2n+k}\ (1/2-\nu-i\tau+\alpha)_n\  (1/2-\nu-i\tau+\beta)_n  \over  n!\ k!\ (1-2(\nu+i\tau))_n}$$

$$\times   {(3\beta-\alpha - 2(\nu+i\tau))_{2n+k}  \over (2(\beta-\alpha - \nu-i\tau))_{2n+k}},$$
where the interchange of the summation and integration can be easily justified by the dominated convergence.  Thus the kernel $ {\cal K}_{\nu+i\tau}^{\alpha,\beta} (x)$ becomes

$$ {\cal K}_{\nu+i\tau}^{\alpha,\beta} (x) = {  x^{ 2(\beta-\alpha - \nu-i\tau)-1} \Gamma(1/2-\beta-\nu-i\tau) \over \pi \Gamma(- 2(\nu+i\tau))  \Gamma(1/2-\alpha+\nu+i\tau) \Gamma( 2(\beta-\alpha - \nu-i\tau))} $$

$$\times \sum_{n,k =0}^\infty  { (- 1)^k x^{2n+k}\ (1/2-\nu-i\tau+\alpha)_n\  (1/2-\nu-i\tau+\beta)_n  \over  n!\ k!\ (1-2(\nu+i\tau))_n}$$

$$\times   {(3\beta-\alpha - 2(\nu+i\tau))_{2n+k}  \over (2(\beta-\alpha - \nu-i\tau))_{2n+k}}+ {  x^{ 2(\beta-\alpha + \nu+i\tau)-1} \Gamma(1/2-\beta+\nu+i\tau) \over \pi \Gamma( 2(\nu+i\tau))  \Gamma(1/2-\alpha-\nu-i\tau) \Gamma( 2(\beta-\alpha + \nu+i\tau))} $$

$$\times \sum_{n,k =0}^\infty  { (- 1)^k x^{2n+k}\ (1/2+\nu+i\tau+\alpha)_n\  (1/2+\nu+i\tau+\beta)_n  (3\beta-\alpha + 2(\nu+i\tau))_{2n+k} \over  n!\ k!\ (1+2(\nu+i\tau))_n\ (2(\beta-\alpha + \nu+i\tau))_{2n+k}}.\eqno(5.7)$$

At the end of this section  we mention some particular cases of the index transformation (1.26) and its inversion (4.15).  Most of them are limit cases which do not satisfy restrictions of Theorem 4 and can be investigated independently.  For instance, when $\nu=0, \alpha=\beta$, it  corresponds to the index transform with the square of the Whittaker function and was considered recently by the author  [11]. Let $\alpha+\beta=0$. Then the fractional integral $I_{0+}^{-\alpha-\beta}$ is the identity operator, and, letting $\beta= -\alpha$,  we have, accordingly,

$$ {\cal K}_{\nu+i\tau}^{\alpha,-\alpha} (x) =   e^{-x} x^{- 4\alpha -1} \bigg[ {  x^{- 2(\nu+i\tau)}\  \Gamma(1/2+\alpha-\nu-i\tau) \over \pi \Gamma(- 2(\nu+i\tau))  \Gamma(1/2-\alpha+\nu+i\tau) \Gamma (- 2(2\alpha+ \nu+i\tau))} \bigg.$$

$$\times    {}_2F_3\left( {1\over 2}-\nu-i\tau+\alpha,\  {1\over 2} -\nu -i\tau-\alpha;\  1- 2(\nu+i\tau),  -2\alpha -\nu-i\tau,\  {1\over 2}-2\alpha-\nu-i\tau;  {x^2\over 4}\right)$$

$$+ {  x^{2(\nu+i\tau) }\    \Gamma(1/2+\alpha+\nu+i\tau) \over \pi \Gamma(2(\nu+i\tau))  \Gamma(1/2-\alpha-\nu-i\tau) \Gamma ( 2(\nu+i\tau-2\alpha))} \  {}_2F_3\left( {1\over 2}+\nu+i\tau+\alpha,\  {1\over 2} +\nu +i\tau-\alpha;\right.$$

$$\left.   1+2(\nu+i\tau), \  \nu+i\tau-2\alpha,\  {1\over 2}+\nu+i\tau-2\alpha;  {x^2\over 4}\right)\bigg],\ x >0.\eqno(5.8)$$
A notable special case is $\alpha= -1/4$,  when the  latter combination of the ${}_2F_3$-hypergeometric functions can be written in terms of the product of the Whittaker's functions $M_{a,b}(z)$.   Indeed, Entries 7.11.1.2, 7.11.1.3,  7.15.1.2 in [2, Vol. III] suggest the following expression for the corresponding kernel

$$ {\cal K}_{\nu+i\tau}^{-1/4, 1/4} (x) =   {e^{-x}\over 2\pi (\nu+i\tau)}\  \bigg[ {  x^{2(\nu+i\tau) }\    \Gamma(1/4+\nu+i\tau) \over  \Gamma^2(2(\nu+i\tau)) \ \Gamma(3/4-\nu-i\tau) } $$

$$\times   {}_1F_1\left( {1\over 4}+\nu+i\tau; \  1+ 2(\nu+i\tau);  x \right)   {}_1F_1\left( {1\over 4}+\nu+i\tau; \  1+ 2(\nu+i\tau);  - x \right)$$

$$- {  x^{- 2(\nu+i\tau)}\  \Gamma(1/4-\nu-i\tau) \over \Gamma^2(- 2(\nu+i\tau))\  \Gamma(3/4+\nu+i\tau) } \bigg.$$

$$\times    {}_1F_1\left( {1\over 4}-\nu-i\tau; \  1- 2(\nu+i\tau);  x \right)   {}_1F_1\left( {1\over 4}-\nu-i\tau; \  1- 2(\nu+i\tau);  - x \right) \bigg]$$

$$=   {e^{-x} x^{-1} \over 2\pi (\nu+i\tau)}\  \bigg[ {   \Gamma(1/4+\nu+i\tau)\  M_{1/4, \nu+i\tau}(x)   M_{-1/4, \nu+i\tau}(x) \over  \Gamma^2(2(\nu+i\tau)) \ \Gamma(3/4-\nu-i\tau) }    $$

$$- {  \Gamma(1/4-\nu-i\tau)\   M_{1/4, -\nu-i\tau}(x)   M_{-1/4, -\nu-i\tau}(x) \over \Gamma^2(- 2(\nu+i\tau))\  \Gamma(3/4+\nu+i\tau) } \bigg],\ x >0.\eqno(5.9)$$
Letting $\nu=0$ in (5.8), we find the kernel ${\cal K}_{i\tau}^{\alpha,-\alpha} (x)$, and the corresponding index transforms modify  the author  results in [8]. 
Finally, for the very particular case $\nu=\alpha=\beta=0$ we arrive at the Lebedev index transform with the square of the Macdonald function [7].  In fact, recalling (1.3), this transform is defined accordingly (see (1.26))

$$(F^0_{0,0} f )(2x)=   e^{2x} \int_{-\infty}^\infty  K_{i\tau}^2 \left(x\right) { f(\tau) \over \cosh(\pi \tau)}\ d\tau.$$
In order to write the corresponding inversion (4.15), let us first calculate the kernel (5.8) with $\alpha=0$.  We have

$$ {\cal K}_{i\tau}^{0,0} (2x) =  { 2 e^{-2x}\over \pi x}\ \cosh(\pi\tau)  \bigg[ {  1 \over  \Gamma^2(- i\tau)  }   {}_1F_2\left( {1\over 2}-i\tau;\  1- 2i\tau,  -i\tau ;  x^2\right) \left({x^2\over 4}\right)^{- i\tau}$$

$$+ {  1 \over  \Gamma^2( i\tau)  }   {}_1F_2\left( {1\over 2}+ i\tau;\  1+ 2i\tau,  i\tau ;  x^2\right) \left({x^2\over 4}\right)^{ i\tau}\bigg].\eqno(5.10)$$
Meanwhile, Entry 7.14.1.4 in [2, Vol. III] presents the expression   of ${}_1F_2 $-functions in terms of the modified Bessel functions of the first kind

$$ {}_1F_2 \left( {1\over 2}  +i\tau ;  \ 1+2i\tau, \   i\tau;\  x^2 \right) =  x \  \Gamma (i\tau)  \Gamma (1+i\tau) \left({t^2\over 4}\right)^{-i\tau} I_{i\tau} \left(x\right) I_{i\tau-1} \left(x\right)$$

$$-  \Gamma^2 (1+i\tau) \left({t^2\over 4}\right)^{-i\tau} I^2_{i\tau} \left(x\right)$$
and,  correspondingly,

$$ {1\over \Gamma^2(i\tau)  }\   {}_1F_2 \left( {1\over 2}  +i\tau ;  \ 1+2i\tau, \   i\tau;\  x^2\right) \ \left({x^2\over 4}\right)^{ i\tau}  =   i\tau \ I_{i\tau} \left(x\right)  \bigg[ x  I_{i\tau-1} \left(x \right)  - i \tau   I_{i\tau} \left(x \right) \bigg].$$
Employing the recurrence relation for the modified Bessel functions [2, Vol. II], the latter formula can be read as follows  

 $$ {1\over \Gamma^2(i\tau)  }\   {}_1F_2 \left( {1\over 2}  +i\tau ;  \ 1+2i\tau, \   i\tau;\  x^2 \right) \ \left({x^2\over 4}\right)^{ i\tau}  =   i\tau x \ I_{i\tau} \left(x\right) \ {d\over dx} \left[ I_{i\tau} \left(x\right)\right] =  i\tau  {x\over 2}  \  {d\over dx} \left[ I^2_{i\tau} \left(x\right)\right].$$
Consequently, the kernel (5.10) becomes

$$ {\cal K}_{i\tau}^{0,0} (2x) =  {  e^{-2x} i\tau \over \pi }\ \cosh(\pi\tau) \ {d\over dx} \bigg[ I^2_{i\tau} \left(x\right) - I^2_{-i\tau} \left(x\right) \bigg],$$
 and  simple functional substitutions suggest the following reciprocal pair of the Lebedev index transforms with the square of the modified Bessel functions 

$$ g(x)=     \int_{-\infty}^\infty  K_{i\tau}^2 \left(x\right)  f(\tau)   d\tau,$$

$$f(\tau)  =   { i \tau\over \pi}  \int_0^\infty  {d\over dt} \bigg[ I^2_{i\tau} \left(x\right) - I^2_{-i\tau} \left(x\right) \bigg]  g(x) dx.$$

\bigskip
\centerline{{\bf Acknowledgments}}
\bigskip

\noindent The work was partially supported by CMUP, which is financed by national funds through FCT (Portugal)  under the project with reference UIDB/00144/2020.   The author sincerely indebted to the referee for a careful reading of the manuscript and valued comments and suggestions which rather improved the presentation of the paper.

\bigskip
\centerline{{\bf References}}
\bigskip
\baselineskip=12pt
\medskip
\begin{enumerate}

 \item[{\bf 1.}\ ] Yu. A. Brychkov, O.I.  Marichev and N.V.  Savischenko, {\it  Handbook of Mellin transforms}. Advances in Applied Mathematics. CRC Press, Boca Raton, FL, 2019.

\item[{\bf 2.}\ ] A.P. Prudnikov, Yu.A. Brychkov and O.I. Marichev, {\it Integrals and Series}. Vol. I: {\it Elementary Functions}, Vol. II: {\it Special Functions}, Gordon and Breach, New York and London, 1986, Vol. III : {\it More special functions},  Gordon and Breach, New York and London,  1990,  Vol. V : {\it Inverse Laplace transforms},  Gordon and Breach, New York and London,  1992.

 \item[{\bf 3.}\ ]  R. Sousa, M.  Guerra and S.  Yakubovich,  {\it Convolution-like structures, differential operators and diffusion processes}. Lecture Notes in Mathematics, 2315. Springer, Cham, 2022.
 
\item[{\bf 4.}\ ]  E.C. Titchmarsh,  {\it Introduction to the theory of Fourier integrals}, Third edition. Chelsea Publishing Co., New York, 1986.

\item[{\bf 5.}\ ]  J. Wimp,  A class of integral transforms,  {\it  Proc. Edinb. Math. Soc.} {\bf 14} (1964), N 2,  33-40. 

\item[{\bf 6.}\ ]  S. Yakubovich and Yu. Luchko, The Hypergeometric Approach to Integral Transforms and Convolutions, {\it Kluwer
Academic Publishers, Mathematics and Applications.} Vol.287, Dordrecht, 1994. 

\item[{\bf 7.}\ ] S. Yakubovich, {\it Index Transforms}, World Scientific Publishing Company, Singapore, New Jersey, London and
Hong Kong, 1996.

\item[{\bf 8.}\ ] S. Yakubovich, Index transforms associated with products of Whittaker's functions. {\it J. Comput. Appl. Math.} Vol. {\bf 148} (2002), no. 2, 419-427.

\item[{\bf 9.}\ ]  S. Yakubovich,  On the Plancherel theorem for the Olevskii transform.  {\it Acta Math Vietnam}.  Vol. {\bf 31} (2006),  249-260.

\item[{\bf 10.}\ ] Nguyen Thanh Hai and  S. Yakubovich (1992). The double Mellin-Barnes type integrals and their applications to convolution theory. {\it World Scientific Intern. Pub!. Series on Soviet and East European Mathematics}. Vol. 6. Singapore, New Jersey, London and Hong Kong, 1992.

\item[{\bf 11.}\ ]  S. Yakubovich,  A new index transform with the square of Whittaker's function. {\it Complex Anal. Oper. Theory}. Vol. {\bf 20} (2026), N 1, 20 pp.

\end{enumerate}

\vspace{5mm}

\noindent S.Yakubovich\\
Department of  Mathematics,\\
Faculty of Sciences,\\
University of Porto,\\
Campo Alegre st., 687\\
4169-007 Porto\\
Portugal\\
E-Mail: syakubov@fc.up.pt\\

\end{document}